  \def\corr{ \vrule height 2pt depth 4pt width 6pt\ }
   \def\obs#1{{\bf (*** #1 ***)} }
 
 \def\corr{}      % Remova esta linha para rodar a versao 1
 \def\obs#1{}     % Remova esta linha para rodar a versao 1

\documentclass[twoside,letterpaper,draft,11pt]{amsart}
\usepackage{amsmath}               % AmSLaTeX
\usepackage{amsthm}                % aggiunge nuovi ambienti tipo teorema

\usepackage{latexsym}

\usepackage{xspace}
\usepackage{amscd}

\title[Associativity of crossed products by partial actions etc.]{Associativity of crossed products by
partial actions, enveloping actions and partial representations}

\author[M.\ Dokuchaev]{M. Dokuchaev}
\address{Departamento de Matem\'atica, Universidade de S\~ao Paulo, Brazil}
\email{dokucha@ime.usp.br}
\author[R.\ Exel]{R. Exel}
\address{Departamento de Matem\'atica, Universidade Federal de Santa Catarina, Brazil}
\email{exel@mtm.ufsc.br}
\urladdr{http://www.mtm.ufsc.br/\~{}exel}
\thanks{{This work was partially supported by CNPq of Brazil
(Proc. 301115/95-8, Proc. 303968/85-0). } }

\newtheorem{teo}{Theorem}[section]
\newtheorem{defi}[teo]{Definition}
\newtheorem{lema}[teo]{Lemma}
\newtheorem{cor}[teo]{Corollary}
\newtheorem{prop}[teo]{Proposition}

\theoremstyle{remark}\newtheorem{remark}[teo]{Remark}
\theoremstyle{remark}\newtheorem{ex}[teo]{Example}

\font\msbm = msbm10
\newcommand{\stimes}{{\hbox {\msbm o}}}

\newcommand{\s}{\sigma}

\newcommand{\D}{{\mathcal D}}

\newcommand{\de}{\delta}

\newcommand{\af}{\alpha}

\newcommand{\bt}{\beta}

\newcommand{\lb}{\lambda}

\newcommand{\gm}{\gamma}

\newcommand{\p}{{\bf Proof. }}

\newcommand{\fim}{\hfill\mbox{$\Box$}}

\newcommand{\A}{{\mathcal A}}
\newcommand{\R}{{\mathcal R}}
\newcommand{\B}{{\mathcal B}}
\newcommand{\I}{{\mathcal I}}
\newcommand{\J}{{\mathcal J}}
\newcommand{\M}{{\mathcal M}}

\newcommand{\X}{{\mathcal X}}
\newcommand{\F}{{\mathcal F}}
\newcommand{\extrassoc}{(L,R)-associative}
\newcommand{\extrassocSubs}{(L,R)-associativity}

\newcommand{\kpar}{K_{\rm par}}     % the partial group algebra
\newcommand{\m}{{}^{-1}}

\def\f{\varphi}
\def\e{\varepsilon}

\newcommand{\G}{\varGamma}

\def\ndv{\ {\mid \kern -0.7 em {\scriptstyle \not}} \ \ }

\def\nd{\ {\mid \kern -0.4 em {\scriptstyle \not}} \ \ }

\begin{document}

\begin{abstract}

Given a partial action $\af$ of a group $G$ on an associative algebra
$\A$ we consider the crossed product $\A\stimes_\af G$.  Using the
algebras of multipliers we generalize a result of  \cite{E0}  on the
associativity of $\A\stimes_\af G$ \corr obtained in the context of $C^*$-algebras. In particular,  we prove that ${\A}
\stimes_{\af} G$ is associative, provided that   $\A$ is semiprime. We also give a criteria for the existence of a global extension of a
given partial action on an algebra and use crossed products to study
relations between partial actions of groups on algebras and
partial representations. As an application we endow partial group
algebras with crossed product structure.

\end{abstract}

\maketitle

\begin{section}{Introduction}\label{sec:intro}

Partial actions of groups appeared  independently in various areas of mathematics, in particular, in the theory of operator algebras
as a powerful tool in their study (see \cite{E-1}, \cite{E0},
\cite{E1}, \cite{Mc}, \cite{QR}). In the most general setting of
partial actions on abstract sets the definition is as follows:

\begin{defi}\label{def1} Let $G$ be a group with identity element $1$ and $\X$ be a set.  A partial action ${\af}$ of
$G$ on  $\X$ is a collection of  subsets ${\D}_g \subseteq \X \; (g \in G)$ and bijections
$ {\af}_g : {\D}_{g\m} \to  {\D}_g$ such that\\

 (i) $\D_1 = \X$ and ${\af}_1$ is the identity map of $\X$;

 (ii) $\D_{(gh)\m} \supseteq {\af}\m_h({\D}_h \cap {\D}_{g\m})$;

 (iii) ${\af}_g \circ {\af}_h (x) = {\af}_{gh}(x)$ for each $x \in {\af}\m_h ({\D}_h \cap {\D}_{g\m})$.\\

\end{defi}

Note that conditions (ii) and (iii) mean that the function ${\af}_{gh}$ is an extension of the function ${\af}_g \circ
{\af}_h $. Moreover, it is easily seen that  (ii) can be replaced by a ``stronger looking'' condition: ${\af}\m_h({\D}_h \cap
{\D}_{g\m}) = \D_{h\m} \cap \D_{h\m g\m}.$ Indeed, it obviously follows from (ii) that ${\af}\m_h({\D}_h \cap
{\D}_{g\m}) \subseteq \D_{h\m} \cap \D_{h\m g\m}.$  Replacing $h$ by $h\m$ and $g$ by $gh$ we have
${\af}\m_{h\m}({\D}_{h\m}
\cap {\D}_{h\m g\m}) \subseteq \D_h \cap \D_{g\m},$ and consequently, $ {\af}_{h\m}({\D}_h
\cap {\D}_{g\m}) \supseteq \D_{h\m} \cap \D_{h\m g\m}.$ By (iii) $\af_{h\m} = \af\m_h$ and we obtain  the desired
equality. Thus, the  conditions (i) - (iii) are equivalent to the following:\\

{\it (i)} $\D_1 = \X$ and ${\af}_1$ is the identity map of $\X${\it ;}

{\it (ii')} ${\af}_g ({\D}_{g\m} \cap {\D}_h) = {\D}_g \cap {\D}_{gh}${\it ;}

{\it (iii')} ${\af}_g({\af}_h(x)) = {\af}_{gh}(x)$, for all $x \in {\D}_{h\m} \cap {\D}_{(gh)\m}.$\\

Let $S(G)$ be the universal semigroup generated by the symbols
$\{ [g] : g \in G \}$ subject to relations:\\

$ a) [g\m][g][h] = [g\m][gh];$

$ b) [g][h][h\m]=[gh][h\m] \; \; \; (g,h \in G);$

$ c) [1]=1,$\\

\noindent where  $1$ denotes also  the identity
element of $S(G)$. Then $S(G)$ is an inverse semigroup  and the partial actions of $G$ on $\X$ are in
one-to-one correspondence with (ordinary) actions of  $S(G)$ on $\X$ \cite{E1}.

In what follows, by an {\it algebra} we shall mean an associative
non-unital (i.e.~non necessarily unital) algebra.  In order to define a partial action
${\af}$ of a group $G$ on a $K$-algebra $\A$ we suppose in
Definition~\ref{def1} that each ${\D}_g \; (g \in G)$ is an ideal of
$\A$ and that every map $ {\af}_g : {\D}_{g\m} \to {\D}_g$ is an
isomorphism of algebras.

Together with the notion of partial actions a generalization  of the concept of crossed product appeared in the theory
of operator algebras (see \cite{E-1}, \cite{E0}, \cite{E2}, \cite{Mc}). For  simplicity  we assume that the
twisting is trivial, so we give the definition  in the  context of  corresponding skew group rings.

\begin{defi} Given a partial action $\af $ of a group $G$ on an
algebra $\A,$ the skew group ring $\A \stimes_{\af} G$ corresponding to
$\af$ is the set of all finite formal sums $ \{ \sum_{g \in G} a_g
\de_g : a_g \in \D_g \},$ where ${\de}_g$ are symbols. Addition is
defined in the obvious way, and multiplication is determined by
$(a_g {\de}_g) \cdot (b_h {\de}_h) = {\af}_g ( {\af}_{g\m}(a_g)b_h)
{\de}_{gh}.$ \end{defi}

Obviously, $\A \ni a \mapsto a \de_1 \in \A \stimes_{\af} G$ is an embedding
which permits us to identify $\A$ with $\A \de_1.$ The first question
which naturally arises is whether or not $\A \stimes_{\af} G$ is
associative.  The associativity of this construction has been proved
in \cite{E-1} in the context of $C^*$-algebras roughly ten years ago.
Since the $C^*$-algebraic proof employed very special properties
of $C^*$-algebras (the existence of approximate units) the
associativity question remained open since then for partial action on
general algebras.

One of our goals is to prove that $\A \stimes_{\af} G$ is always
associative if  $\A$ is semiprime (Corollary~\ref{idempotent2}).  
\corr Recall  that a unital algebra $\A$ is called semiprime  
if $\A$ has no non-zero nilpotent ideal 
(for other equivalent definitions  see \cite[Theorem 2.6.17]{Rowen}).  
\corr More generally, given a partial action of $G$ on a unital algebra $\A,$ the skew group ring   $\A \stimes_{\af} G$ is associative provided that each ideal $\D_g$   is idempotent or non-degenerate \corr (Corollary~\ref{idempotent1.5}).   We say that an ideal $\I$ of  $\A$ is non-degenerate  if for every non-zero element $a\in \I$ there exists  $b\in\I$ such that either $ab\neq0$ or $ba\neq0$. \corr It is easily proved that a  unital algebra $\A$ is semiprime if and only if  each non-zero ideal of $\A$ is  non-degenerate  (see Proposition~\ref{semiprime}).

 Note that in $C^*$-algebraic context the ideals $\D_g$ are supposed to be 
closed by definition, and it is known that each closed ideal in a
$C^*$-algebra is an idempotent ideal (see \cite[Theorem
V.9.2]{FellDoran}).

We also show that in general ${\A}
\stimes_{\af} G$ is not associative (see Section~\ref{sec:associativity}).
In Section~\ref{sec:envelope} we give a  complete answer for the problem
of the existence of a ``global extension'' (called enveloping action)
for a given partial action on a unital algebra. More precisely, an
enveloping action
 exists if and only if each ideal $\D_g$ is generated by a central idempotent 
(see Theorem~\ref{envelopeexistence}). We
introduce the notion of enveloping action with certain minimality condition in 
order to guarantee its uniqueness up to
equivalence (see Section~\ref{sec:envelope} for the definitions).
  We also prove
 that if $(\B,\bt)$ is an enveloping action for a
partial action $(\A,\af)$ such that  $\A$ and $\B$ are unital algebras
then the corresponding skew group rings,
namely $\A\stimes_\af G$ and $\B\stimes_\bt G$, are Morita equivalent.

 Another relevant concept which appeared in the theory of operator
algebras is the notion of partial representation of a group on a
Hilbert space, introduced independently in \cite{E1} and
\cite{QR}. Partial representations underlie important algebras
generated by partial isometries.  Among the most interesting cases are
the Cuntz-Krieger algebras \cite{CK} investigated in \cite{E2}
and \cite{EL} from the point of view of partial representations. The
purely algebraic study of partial representations began in \cite{DEP}
and was continued in \cite{DZh0} (see also \cite{DZh}). Similarly to
the case of (usual) representations of groups, there exists an algebra
$\kpar (G)$, called the partial group algebra of $G$ over $K$, which
governs the partial representations of $G$ (see \cite{DEP}). It is
exactly the semigroup algebra $KS(G)$.  It turns out that $\kpar(G)$
keeps much more information about the structure of $G$ than does $KG$
(see \cite{DEP}) and this makes them especially interesting with respect
to the classical isomorphism problem intensively investigated for group rings 
(see \cite{Sehgal}).

In the final Section~\ref{sec:reprcross} we use crossed products  to relate partial actions with partial representations.
This leads to a one-to-one correspondence (see Theorem~\ref{correspondence1}) when working with so-called ``elementary''
partial representations, which are structural blocks of the irreducible partial representations (see \cite{DZh0}).
This permits to look at matrix algebras as crossed products  by
partial actions (see Corollary~\ref{isoBIS}),
  covering this way a
substantial part of elementary gradings (see Corollary~\ref{gradings}). Moreover, starting with the partial representation $G
\ni g \mapsto [g] \in \kpar (G),$ we  endow
$\kpar (G)$  with crossed product structure  (see Theorem~\ref{pgr}).

\end{section}

\begin{section}{The algebra of multipliers}\label{sec:multipliers}

Let $K$ be a field, $\A$  an associative  $K$-algebra with unity element  and $\I$   an ideal of $\A$.
Take an element $x \in \A$ and consider the left and right multiplications of   $\I$ by $x$:  $L_x: \I \ni a
\mapsto xa \in  \I, \; \;\; R_x:  \I \ni a \mapsto ax \in \I.$ Then $L = L_x$
and
$R = R_x$ are linear transformations of $\I$ such that  the following properties are satisfied for all $a, b \in \I$:\\

\noindent (i)   $L(ab) = L(a)b$;\\
(ii)  $R(ab) = aR(b)$;\\
(iii) $R(a)b = a L(b)$.\\

These properties are obvious consequences of the associativity of $\A$.\\

\begin{defi}\label{multipalgebra} The algebra of multipliers (see e.g. \cite[3.12.2]{Fillmore}) of an algebra  $\I$  is
the set $M(\I)$ of all ordered pairs $(L,R)$, where $L$ and
$R$ are linear transformations of $\I$ which satisfy the properties (i) - (iii). For $(L,R), (L',R') \in M(\I)$ and
$\af \in K$
the operations are given by $\af (L,R) = (\af L, \af R), \; (L,R) + (L',R') = (L+L', R+R'), \; (L,R) (L', R') =
(L \circ L', R'
\circ R).$ We say that $L$ is a left multiplier and $R$ is a right multiplier of $\I$.\\
\end{defi}

 It is immediately verified that $M(\I)$ is an associative algebra with unity element $(L_1, R_1)$, where $L_1$ and $R_1$ are
identity maps (which in the case of an ideal $\I$ in a unital algebra $\A$ can be considered as multiplications by the unity
element of $\A$ from left and right, respectively). Define the map $\phi: \I \to M(\I)$  by putting $\phi(x) = (L_x,R_x),
\; x \in \I$. This is a homomorphism of algebras since it is $K$-linear and, moreover, $L_{xy} = L_x \circ L_y, R_{xy} = R_y
\circ R_x$, which gives
$\phi(xy) =  (L_x \circ L_y, R_y \circ R_x ) = \phi (x) \phi (y).$ 

\begin{defi} We shall say that an algebra $\I$ is non-degenerate if
the map $\phi: \I \to M(\I)$ mentioned above is injective.
\end{defi}

In general the kernel of $\phi$ is the intersection of the left
annihilator of $\I$ in $\I$ with its right annihilator in $\I$.
Therefore $\I$ is non-degenerate if and only if for every non-zero
element $a\in \I$ there exists $b\in\I$ such that either $ab\neq0$ or
$ba\neq0$.

More generally, if $\I$ is an ideal in an algebra $\A$ then one may consider the homomorphism $\psi: \A \ni a \mapsto   (L_a,R_a) \in M(\I),$
whose kernel  is the intersection of the left annihilator of $\I$ in $\A$ with its right annihilator in
$\A$.

%\newpage

\begin{prop}\label{M(I)} The following statements hold:

\noindent (i) $\phi(\I)$ is an ideal of $M(\I).$

\noindent (ii) $\phi: \I \to \M(\I)$ is an isomorphism if and only if $\I$ is a unital algebra.

\end{prop}

\p (i) Take $x \in \I$ and let $ (L,R)$ be an arbitrary element of $M(\I)$. Then $(L_x, R_x) (L,R) = (L_x \circ L, R \circ
R_x)$ and for $a \in \I$ we have $L_x  (L(a)) = x L(a) = R(x)a = L_{R(x)}(a)$. Moreover, $R (R_x (a)) = R(ax) = a
R(x) = R_{R(x)}(a)$. Hence, $(L_x, R_x) (L, R) = ( L_{R(x)}, R_{R(x)} ) \in \phi(\I)$, as $R(x) \in \I$. Similarly, $(L,R)
(L_x,R_x) = (L_{L(x)}, R_{L(x)}) \in \phi(\I)$.\\

\noindent (ii) The ``only if''  part is trivial. For the ``if'' part we have that $\phi(1) \in \phi(\I)$ is the identity
element of $M(\I)$ and, consequently, $\phi(\I) = M(\I)$. Obviously, $\phi$ is injective in this case and thus
$\I \cong M(\I)$.
\fim\\

Let $\I$ be any (preferably non-unital) algebra.  Given $(L,R)$ and
$(L',R')$ in $\M(\I)$ we shall be concerned with the validity of the
formula 
\begin{equation}\label{extraassoc}
R' \circ L = L \circ R'.
\end{equation}

If $x$ and $x'$ belong to an algebra which contains $\I$ as an ideal
and $(L,R)=(L_x,R_x)$ and $(L',R')=(L_{x'},R_{x'})$, this formula
will clearly hold as a consequence of associativity.  However this is
not always the case:  for a very drastic counter-example one could take  $\I$
to be any vector space equipped with the trivial multiplication
  $$
  xy \equiv 0,\qquad \forall x,y\in\I.
  $$
  Any pair $(L,R)$ of linear operators on $\I$ would constitute a
multiplier of $\I$ and one would clearly not expect (\ref{extraassoc})
to hold!

\begin{defi}\label{extraassocdef} An algebra $\I$  is said to be {\extrassoc} if, given any two
multipliers $(L,R)$ and $(L',R')$ in $\M(\I)$, one has that $R' \circ
L = L \circ R'.$ \end{defi}

The following result lists two sufficient conditions for {\extrassocSubs}.

\begin{prop}\label{condforextraassoc} The algebra $\I$ is
{\extrassoc} whenever any one of the following conditions hold:

\noindent (i) $\I$ is non-degenerate, or

\noindent (ii) $\I$ is idempotent.
\end{prop}

\p Let  $(L,R),(L',R')\in\M(\I)$.  Given  $a,b\in\I$ we have that
  $$
  R(L'(a))b =
  L'(a)L(b) =
  L'(aL(b)) =
  L'(R(a)b) =
  L'(R(a))b.
  $$
  This shows that $R(L'(a))-L'(R(a))$ lies in the left annihilator of
$\I$.  By a similar calculation one shows that $R(L'(a))-L'(R(a))$ lies
in the right annihilator of $\I$ as well.  Therefore, under assumption
$(i)$, one must have that $R(L'(a))=L'(R(a))$, for all $a$ in $\I$.

Next suppose  we are given $a_1,a_2\in\I$.  Letting
$a=a_1a_2,$  notice that
  $$
  R(L'(a)) =
  R(L'(a_1a_2)) =
  R(L'(a_1)a_2) =
  L'(a_1)R(a_2) $$$$
  = L'(a_1R(a_2)) =
  L'(R(a_1a_2)) =
  L'(R(a)).
  $$
  Assuming $(ii)$ we have that  every element of $\I$ is a sum of terms of the
  form $a_1a_2$, whence the conclusion.
\fim\\

In the next section we will be considering partial actions such that
all of the ideals $D_g$ are assumed to be {\extrassoc}.  It is
therefore useful to have a means of deciding when all ideals of a
certain algebra possess this property.  The following  easy result goes
in that direction.

\corr \begin{prop}\label{semiprime} Let $\A$ be a unital algebra.  The following are
equivalent:

\begin{itemize}
\item[\it (i)] Every non-zero  ideal of $\A$ is non-degenerate,
\item[\it (ii)] Every non-zero  ideal of $\A$ is either  idempotent or non-degenerate,
 
\item[\it (iii)] Every non-zero ideal of $\A$ is right-non-degenerate
(we say that  $\I$ is right-non-degenerate if for every non-zero element
$a\in\I$ one has that $a\I\neq\{0\}$),
 
\item[\it (iv)] Every non-zero ideal of $\A$ is left-non-degenerate
(defined in a similar way  with $\I a\neq\{0\}$),
 
\item[\it (v)] $\A$  is semiprime.
\end{itemize}

\medskip \noindent In this case every ideal of $\A$ is {\extrassoc}.
\end{prop}

\def\starArrow{\mathrel{{\Rightarrow\kern-12pt {}^{{}^{\textstyle*}}}}}

\p  \corr In  the following implications the
ones marked with a star are self-evident:
  $$(i) \starArrow (ii) \Rightarrow  (v) \Rightarrow (iii) \starArrow (i).$$
  Moreover $(iv)$ can be substituted for $(iii)$, by symmetry, so we
really only need to worry about \corr  $(ii) \Rightarrow (v) \Rightarrow (iii).$

\corr If $\A$ has a non-zero nilpotent ideal then  there is a non-zero ideal $\I$  in $\A$ whose square is zero. Then $\I$ is neither idempotent nor non-degenerate which shows  $(ii) \Rightarrow (v).$

Assuming $(v)$, and  arguing by contradiction, let $\I$ be an
ideal possessing a non-zero element $a$ such that $a\I=\{0\}$.  Then
the ideal generated by $a$, namely  $\J=\A a\A$, is non-zero
since we are assuming $\A$ to be unital.  However
  $$
  \J^2 = \A a(\A \A a)\A \subseteq \A a \I\A = \{0\},
  $$
  violating $(v).$
\fim\\

Let  $\pi : \I \to \J$ be an isomorphism of $K$-algebras. Then it is easy to see that for $(L,R) \in M(\I)$ the
pair $(\pi \circ L \circ {\pi}^{-1}, \pi \circ R \circ {\pi}^{-1})$ is an element of $M(\J)$ and we obviously have the
following:

\begin{prop}\label{iso} The map  $\bar{\pi}: M(\I) \to M(\J),$ defined by  $$ \bar{\pi} (L,R) = (\pi \circ L \circ
{\pi}^{-1}, \pi \circ R \circ {\pi}^{-1}),$$ is an isomorphism of $K$-algebras.\\
\end{prop}

% Indeed, for $a, b \in \J$ we have
% that

% $$\pi \circ L \circ {\pi}^{-1} (ab) = \pi  (L[\pi^{-1}(a) \pi^{-1}(b)]) = \pi(L[\pi^{-1}(a)] \pi ^{-1}(b)) =
% (\pi \circ L \circ \pi^{-1}) (a)  b,$$

% \noindent and

% $$\pi \circ R \circ {\pi}^{-1} (ab) = \pi (R[\pi^{-1}(a) \pi^{-1}(b)]) = \pi(\pi^{-1}(a) R[\pi ^{-1}(b))] =
% a (\pi \circ R \circ \pi^{-1}) (b).$$

% \noindent Moreover,

% $$  (\pi \circ R \circ \pi^{-1}) (a) b = \pi( R [\pi^{-1}(a)] \pi^{-1}(b)) = \pi( \pi^{-1}(a) L[\pi^{-1}(b)]) =
% a (\pi \circ L \circ \pi^{-1}) (b).$$

% \noindent Hence,  $(\pi \circ L \circ {\pi}^{-1}, \pi \circ R \circ
% {\pi}^{-1})$  satisfies properties (i)-(iii) of the definition of a multiplier of $\J$.

\end{section}

\begin{section}{The associativity question}\label{sec:associativity}

We are now ready to present an answer to  the
associativity question:\\

\begin{teo}\label{idempotent1} If $\A$ is an algebra and $\af $ is a
partial action of a group $G$ on $\A$ such that each $\D_g$ $(g \in
G)$ is {\extrassoc}, then the skew group ring $\A\stimes_\af G$
is associative.  \end{teo}

\p Obviously, $\A\stimes_\af G$ is associative if and only if

\begin{equation}\label{hgf}
(a \delta_h b \delta_g ) c \delta_f = a \delta_h ( b \delta_g c \delta_f)
\end{equation}

\noindent for arbitrary $h, g, f \in G$ and $a \in \D_h, b \in \D_g, c \in \D_f$. We compute first the left
hand side of the above equality. We have

$$ (a \delta_h b \delta_g ) c \delta_f =  \af _h ( \af_{h^{-1}}(a) b)  \delta_{hg} c  \delta_f = \af _{hg} \{
\af^{-1}_{hg}  [\af _h ( \af_{h^{-1}}(a) b) ] c \}     \delta_{hgf}.$$

\noindent We see that $\af_{h^{-1}}(a) b \in \D_{h^{-1}} \cap \D_g$ implies $\af_h (\af_{h^{-1}}(a) b) \in
\af_h( \D_{h^{-1}} \cap \D_g) = \D_{h} \cap \D_{hg}.$ It follows that $\af^{-1}_{hg}  [\af _h (
\af_{h^{-1}}(a) b) ] = \af_{g^{-1}}( \af_{h^{-1}} [\af _h ( \af_{h^{-1}}(a) b) ]) =  \af_{g^{-1}}(
\af_{h^{-1}}(a) b).$ Since  this element belongs to $\D_{g\m} \cap \D_{g\m h\m}$, we can also split $\af_{hg}$, which
gives

$$ (a \delta_h b \delta_g ) c \delta_f = \af_h[ \af_g \{ \af_{g^{-1}}( \af_{h^{-1}}(a) b) c\} ] \delta_{hgf}.$$

\noindent Comparing with

$$ a \delta_h (b \delta_g  c \delta_f ) =   a \delta_h \af _g ( \af_{g^{-1}}(b) c)  \delta_{gf} = \af _{h} [
\af_{h^{-1}} (a) \af _g ( \af_{g^{-1}}(b) c) ]     \delta_{hgf},$$

\noindent and applying $\af_{h^{-1}}$, we obtain that (\ref{hgf}) holds if and only if

$$ \af_g \{ \af_{g^{-1}}( \af_{h^{-1}}(a) b) c\} = \af_{h^{-1}} (a) \af _g ( \af_{g^{-1}}(b) c) $$

\noindent is verified for all $ a \in \D_h, b \in \D_g, c \in \D_f$. Because $\af_{h^{-1}} : \D_h
\longrightarrow \D_{h^{-1}}$ is an isomorphism, $\af_{h^{-1}}(a)$ runs over $\D_{h^{-1}}$ and, consequently, the above
condition is equivalent to the following:

\begin{equation}\label{X}
\af_g \{ \af_{g^{-1}}( a b) c\} = a \af _g ( \af_{g^{-1}}(b) c)
\end{equation}

\noindent  for every $ a \in \D_{h^{-1}}, b \in \D_g, c \in \D_f.$ If $h = f = 1$, then
$\D_h   = \D_f = \A$ and thus $\A\stimes_\af G$ is associative if and only if (\ref{X}) holds for
 arbitrary $g \in G,  a, c \in \A$ and $   b \in \D_g.$  It is equivalent to say that

\begin{equation}\label{XX}
(\af_g \circ R_c \circ \af_{g^{-1}}) \circ L_a = L_a \circ (\af_g \circ R_c \circ \af_{g^{-1}})
\end{equation}

\noindent is valid on $\D_g$ for every $g \in G$ and all  $a, c \in \A$.

Consider $R_c$ as a right multiplier of $\D_{g^{-1}}$ and $L_a$ as a
left multiplier of $\D_g$.  By Proposition~\ref{iso} we have that
$\af_g \circ R_c \circ \af_{g^{-1}}$ is a right multiplier of
$\D_g$. Hence if $\D_g$ is {\extrassoc} then
(\ref{XX}) holds. \fim\\

\begin{cor}\label{idempotent1.5} If  $\af $ is a partial action of a group $G$ on an algebra $\A$ such that each $\D_g$ $(g \in
G)$ is either idempotent or non-degenerate, then the skew group ring $\A\stimes_\af G$
is associative.  
\end{cor}

\p Directly follows from   Proposition~\ref{condforextraassoc} and Theorem~\ref{idempotent1}. \fim\\

For ease of reference it would be useful to introduce the following
terminology:

\begin{defi} We say that an algebra $\A$ is {\it strongly associative} if for any  group $G$ and an
arbitrary partial action $\af$ of $G$ on $\A$ the skew group ring $\A \stimes_{\af} G$ is associative.
\end{defi}

 As an  immediate consequence of Proposition~\ref{semiprime} and Corollary~\ref{idempotent1.5} we have:

\begin{cor}\label{idempotent2}
 A semiprime algebra is strongly associative.  
\end{cor}

 The following is an easy example which shows that  $\A\stimes_\af G$ is not associative in general.\\

\begin{ex}\label{counter}

Let $\A$ be a four-dimensional $K$-vector space with basis $\{ 1, t, u, v \}$. Define the multiplication on
$\A$ by setting $u^2 = v^2 = uv = vu = t u = u t = t^2 = 0, t v = v t = u$ and $1 a = a 1 = a$ for each $a \in \A$.
Then $\A$ is an associative $K$-algebra with unity. Let $G = \langle g : g^2 = 1 \rangle$ and $\I$ be the ideal of $\A$
generated by $v$ (it is the subspace generated by $u$ and $v$). Consider
the partial  action $\af$ of $G$ on $\A$ given by $\D_g = \I$, $\af_g : u \mapsto v, v \mapsto u$ (by  definition
$\D_1 = \A$ and $\af_1$ is the identity map of $\A$). Then the skew group ring $\A \stimes_{\af} G$ is not associative.
More precisely, taking $x = t {\de}_1 + u {\de}_g$ we have that $(xx)x = 0$ and $x(xx) = u {\de}_g$, so that $\A \stimes_{\af} G$ does   not
even have associative powers.\\
\end{ex}

  \corr An important class of non-semiprime algebras is formed by the group algebras
$K G$ of finite groups $G$ with $char K $ dividing the order of $G$ (for more general information see \cite[Theorems 2.12 and 2.13]{Passman}). It is  easy to check that if $char K =2$ then the algebra of the above example is isomorphic to the group algebra of the Klein four-group over $K$. On the other hand, it can be verified
that the group algebra of the cyclic group of order $2$ over a field of characteristic $2$ is strongly associative.  Thus it
would be interesting to characterize the strongly associative group algebras. Another classical example   of \corr  non-semiprime
algebras is given by the algebra $T(n,K)$ of upper-triangular $n\times n$-matrices over $K$.

\begin{prop}\label{triangular} The algebra $\A = T(n,K)$ is strongly associative if and only if $n \leq 2.$
\end{prop}

\p If $n = 1$ then $\A = K$ is obviously strongly associative. If $n
=2,$ then the only non-idempotent ideal of $\A$ is its Jacobson
radical $\R(\A)$, which is one-dimensional over $K$. So all
multipliers of $\R(\A)$ commute and hence it is {\extrassoc}.

Suppose that $n \geq 3$ and let $G = \langle g \rangle$ be the infinite cyclic group. Denote by $e_{i,j}$  the
elementary $n \times n$-matrix whose unique non-zero entry equals $1$ and is placed at the intersection of the $i$-th
row and
$j$-th column. Take
$\D_{g\m} = e_{1,n-1}K \oplus e_{1,n}K,$ $\D_g = e_{1,n}K \oplus
e_{2,n}K$ and $\D_{g^m} = e_{1,n}K$
for each $m$ with $|m| \geq 2.$
Define $\af_g : \D_{g\m} \to \D_g$ by $\af_g (x e_{1,n-1} + y e_{1,n})
= y e_{1,n}+ x e_{2,n}$ $(x,y \in K)$ and for $m \geq 2$ let
$\af_{g^m} : \D_{g^{-m}}= \D_{g^m} \to \D_{g^m}$ be the identity
map. An easy verification shows that we have defined a partial action
$\af$ of $G$ on $\A$.  We see that

\begin{align*} & (e_{1,1} \de_1 \cdot e_{2,n} \de_g) e_{n-1,n} \de_1 = (e_{1,1} \cdot e_{2,n} \de_g) e_{n-1,n} \de_1 = 0,
\end{align*}

\noindent as $e_{1,1} \cdot e_{2,n} = 0.$ On the other hand

\begin{align*} & e_{1,1} \de_1  ( e_{2,n} \de_g \cdot e_{n-1,n} \de_1)  = e_{1,1} \de_1 (\af_g ( \af_{g\m}(e_{2,n})e_{n-1,n})
\de_g)   = \\ &=  e_{1,1} \de_1 ( \af_g ( e_{1,n-1} \cdot e_{n-1,n}) \de_g) = e_{1,1} \de_1 \cdot \af_g (e_{1,n}) \de_g =
e_{1,1} \de_1 e_{1,n} \de_g = \\&= (e_{1,1} \cdot e_{1,n}) \de_g = e_{1,n} \de_g \neq 0,
\end{align*}

\noindent so that $\A \stimes_{\af} G$ is non-associative. \fim

\end{section}

\begin{section}{Enveloping actions}\label{sec:envelope}

Natural examples of partial actions can be obtained by restricting  (global) actions to non-necessarily invariant
subsets (ideals in our case). More precisely, suppose that a group $G$ acts on an algebra $\B$ by automorphisms
$\bt_g : \B \to \B$  and  let $\A$ be an ideal  of $\B.$  Set $\D_g = \A \cap \bt_g(\A)$ and let $\af_g$ be the restriction of
$\bt_g$ to $\D_{g\m}.$   Then it is easily verified that we have a partial action $\af = \{ {\af}_g : {\D}_{g\m} \to
{\D}_g : g \in G \}$ of $G$ on $\A$. We shall say that
$\af$ is a {\it restriction} of $\bt$ to $\A.$ We want to know  circumstances under which a given partial action can be
obtained ``up to equivalence'' as the restriction of  a (global) action. If $\B_1$ is the subalgebra of $\B$ generated by
$\cup_{g \in G} \bt_g(\A),$ it may happen that $\B \neq \B_1$  and, since $\B_1$ is invariant with respect to $\bt,$ we see
that $\af$ can be obtained as a restriction of  an action of $G$ to
$\B_1$ which is a  proper subalgebra of  $\B.$ Thus it is reasonable to require that $\B = \B_1$ and in this case we say that $\af$ is an {\it admissible
restriction} of $\bt.$ The notion of the {\it equivalence} of partial
actions  is  defined as follows:

\begin{defi}\label{equiv} We say that a partial action $\af = \{ {\af}_g : {\D}_{g\m} \to  {\D}_g : g \in G \}$  of a
group $G$ on an algebra $\A$ is equivalent to the partial action $\af' = \{ {\af'}_g : {\D'}_{g\m} \to  {\D'}_g : g \in G
\}$ of $G$ on an algebra $\A'$ if there exists an algebra isomorphism
$\f: \A \to \A'$ such that for every $g \in G$ the following two conditions hold:

{\it (i)} $\f (\D_g) = \D'_g${\it ;}

{\it (ii)}  ${\af'}_g \circ \f (x) = \f \circ {\af}_g (x)$ for  all $x \in \D_{g\m}.$\\

\end{defi}

We shall deal with {\it enveloping actions}:

\begin{defi}\label{envelopeaction} An action $\bt$ of a group
$G$ on an algebra $\B$ is said to be an  enveloping action for the  partial action $\af$ of $G$ on an algebra $\A$ if $\af$ is
equivalent to an admissible restriction of $\bt$ to an ideal of $\B.$
\end{defi}

In other words,  $\bt$ is an enveloping action for $\af$ if there
exists an algebra  isomorphism $\f $ of $\A$  onto an ideal  of  $\B$ such that   for all $g \in G$ the
following three properties  are satisfied:

{\it (i')}  $\f(\D_g) = \f(\A) \cap \bt_g(\f(\A))${\it ;}

{\it (ii')} $ \f \circ \af_g  (x) = \bt_g \circ \f  (x)$ for each $x \in   \D_{g\m},$

{\it (iii')} $\B$ is generated by $\cup_{g \in G} \bt_g(\f(\A)).$\\

Thus a general problem is to decide whether or not a given partial action possesses an enveloping action.

With respect to  the results of Section~\ref{sec:associativity}   we observe the following:

\begin{prop}\label{envelopecross}

 If $\bt$ is an  action of a group $G$ on an algebra $\B,$ which is enveloping for the partial action $\af$ of $G$
on an algebra $\A,$ then the skew group ring  $\A \stimes_{\af} G$ has an embedding into $\B \stimes_{\bt} G$. In particular,
$\A \stimes_{\af} G$ is associative.
\end{prop}

\p Obvious. \fim

Thus it follows from Example~\ref{counter} (or Proposition~\ref{triangular})
that not every partial action admits an enveloping action. We need the following easy fact:

\begin{lema}\label{unity} Suppose that $\A$ is an algebra which is a
(non-necessarily direct) 
 sum of a finite number of ideals, each of
which is a unital algebra. Then $\A$ is a unital algebra.
\end{lema}

\p By induction on the number of summands it is enough to consider the case with two ideals: $\A = \I +\J$. Let
$1_{\I}$ and $1_{\J}$ be the unity elements of  $\I$ and $\J$ respectively. Then $1_{\I}$ and $1_{\J}$ are central
idempotents of $\A$ and  $1_{\I} + 1_{\J} - 1_{\I} \cdot 1_{\J}$ is the unity of $\A.$ \fim

The following result is inspired  by F.~Abadie's thesis
\cite{Abadie} (see also \cite{AbadieTwo}):

\begin{teo}\label{envelopeexistence}
 Let $\A$ be a unital algebra. Then a partial action   $\af $   of a group $G$  on     $\A$ admits an enveloping action $\bt$
if and only if   each ideal $\D_g$ $(g \in G)$ is a unital algebra. Moreover, $\bt,$ if it exists, is unique up to
equivalence.
\end{teo}

\p The ``only if'' part is trivial, because if $\bt$ exists and $\f : \A \to \B$ is the monomorphism giving the equivalence,
then  $\f (\D_g) = \f(\A) \cap \bt_g(\f(\A))$ is clearly a unital algebra for each $g \in G.$

 For the ``if'' part we suppose that each ideal  $\D_g$ $(g \in G)$ is a unital algebra. It means that for every $g \in
G$ there exists a central idempotent $1_g$ of $\A$ such that $\D_g = 1_g \A.$

Let $\F = \F(G,\A)$ be the Cartesian product of the copies of $\A$ indexed by the elements of $G,$ that is, the algebra of all
functions  of $G$ into $\A.$ For convenience of notation $f(g)$ will be also written as  $f|_g$  $(f \in \F, g \in G).$

For $g \in G$ and $f \in \F$ define $\bt_g(f) \in \F$ by the formula:
$$\bt_g (f)|_h = f(g\m h), \;h \in G.$$ It is easily verified that $f
\mapsto \bt_g(f)$ defines an automorphism $\bt_g$ of $\F$ and hence
$\bt = \{\bt_g: \F \to \F: g \in G\} $ is an action 
of $G$ on $\F$.

 It
is easy to see that for each $g, h \in G$ the idempotent $1_g 1_h$ is the
unity element of the algebra $\D_g \cap \D_h,$ which means that $\D_g \cap \D_h = 1_g 1_h \A$. Because $\af$ is
a partial action,  the equality
${\af}_g ({\D}_{g\m} \cap {\D}_h) = {\D}_g \cap {\D}_{gh}$ obviously implies

\begin{equation}\label{triviality}
\af_g (1_{g\m} 1_h) = 1_g 1_{gh}.
 \end{equation}

For any $a \in \A$ the element $a 1_g$  belongs to $\D_g$ and the formula $$\f(a)|_g = \af_{g\m}(a 1_g), \; g \in G,$$ defines
a monomorphism $\f : \A \to \F.$

 Let $\B$ be the subalgebra of $\F$ generated by  $\cup_{g \in G} \bt_g(\f(\A)),\; (g \in G).$
Our purpose  is to show  that the restriction of $\bt$ to $\B$   is an enveloping action for $\af$. We denote this
restriction by the same symbol $\bt.$  We start by checking property (ii') of  Definition~\ref{envelopeaction}.

For $g,h \in G$ and $a \in \D_{g^{-1}}$ we have $\bt_g(\f(a))|_h = \f (a)|_{g\m h} = \af_{h\m g} (a 1_{g\m h})$ and
$\f(\af_g(a))|_h = \af_{h\m}(\af_g(a) 1_h).$ Thus (ii') is satisfied if and only if the following equality  holds for all
$g,h \in G$ and every $a \in \D_{g\m}$:

\begin{equation}\label{alfa}
\af_{h\m g}(a 1_{g\m h}) = \af_{h\m}(\af_g (a) 1_h).
\end{equation}

\noindent Observing that $a \cdot 1_{g\m h} \in \D_{g\m} \cap \D_{g\m h}$ we can split $\af_{h\m g}$ in the left hand side and
using (\ref{triviality}) we obtain:

\begin{align*} & \af_{h\m g}(a 1_{g\m h}) = \af_{h\m} ( \af_g(a 1_{g\m h})) =  \af_{h\m} ( \af_g(a 1_{g\m} 1_{g\m h})) =  \\&=
 \af_{h\m} ( \af_g(a) \af_g(1_{g\m} 1_{g\m h})) =   \af_{h\m} ( \af_g(a) 1_g 1_h) =  \af_{h\m} ( \af_g(a)  1_h),
\end{align*}

\noindent as $1_g$ is the unity of $\D_g.$

Next we show that

\begin{equation}\label{deltas}
 \f(\D_g) = \f(\A) \cap \bt_g(\f(\A)),
\end{equation}

\noindent for all $g \in G.$ An element from the right hand side can be written as $\f(a) = \bt_g(\f(b))$ for some $a,b \in
\A.$ Then for each $h \in G$ the equality $\f(a)|_h = \bt_g(\f(b))|_h $ means that

\begin{equation}\label{triviality2}  \af_{h\m}(a 1_h) =   \f(b)|_{g\m h} = \af_{h\m g}(b 1_{g\m h}).
\end{equation}

\noindent Taking $h=1$ this gives $a = \af_g(b 1_{g\m}) \in \D_g$ and, consequently, $\f(\D_g) \supseteq \f(\A) \cap
\bt_g(\f(\A)).$  For the reverse inclusion, given an element $a \in \D_g$ we need to find $b \in \A$ so that
(\ref{triviality2}) holds. For $b = \af_{g\m}(a)$ the right hand side of (\ref{triviality2}) is $ \af_{h\m g}(\af_{g\m}(a)
1_{g\m h})$ which is equal to the left hand side of  (\ref{triviality2}) in view of   (\ref{alfa}). Hence
(\ref{deltas}) follows and condition (i') is also satisfied.

 In order to show that $\bt$ is an  enveloping action for $\af$ it remains to prove that $\f(\A)$ is an ideal of $\B.$
To see this it is enough to check that $\bt_g(\f(a)) \cdot \f(b), \f(b)  \cdot \bt_g(\f(a)) \in \f(\A)$ for all $g \in G$ and
$ a,b \in \A.$  For $h \in G,$  using (\ref{alfa}), we have

\begin{align*} & \bt_g(\f(a))|_h \cdot \f(b)|_h = \f(a)|_{g\m h} \cdot \f(b)|_h = \af_{h\m g} (a 1_{g\m h} ) \cdot
\af _{h\m}(b 1_h) = \\ & = \af_{h\m}(\af_g(a 1_{g\m})1_h) \cdot \af_{h\m}(b 1_h) = \af_{h\m} ( \af_g(a 1_{g\m} ) b 1_h) = \f(\af_g(a 1_{g\m} )b)|_h.
\end{align*}

\noindent  Thus $\bt_g(\f(a)) \cdot \f(b) = \f(\af_g(a 1_{g\m} )b) \in \f(\A)$ and  similarly   $ \f(b) \cdot \bt_g(\f(a))  =\\
\f(b \af_g(a 1_{g\m} )) \in \f(\A),$ as desired.

We shall prove now the uniqueness of the  enveloping action. Suppose that  $\bt'$ is another   action of $G$ on
an algebra $\B'$, which is   enveloping for $\af$. Let  $\f'$ be the corresponding
embedding of
$\A$ into  $\B'.$ The admissibility incorporated in Definition~\ref{envelopeaction}  means that $\B'$ is the sum of
the ideals $\bt'_g(\f'(\A)),\;g\in G.$  Thus an element of $\B'$ can be written as a finite sum $   \sum_i
\bt'_{g_i} (\f'(a_i))$ with $ g_{i} \in G$ and $a_i \in \A.$  Define a map $\phi: \B' \to \B$ by
$\bt'_g(\f'(a)) \mapsto \bt_g(\f(a)), \;g\in G, a\in \A.$ We have to show, of course, that $\phi$ is well
defined. Suppose that $   \sum_{i=1}^s \bt'_{g_i} (\f'(a_i)) =0.$ We want to be sure that $\sum_{i=1}^s \bt_{g_i} (\f(a_i))
= 0.$

 For all $h \in G$ and $ a \in \A$ we have
$\sum_i \bt'_{g_i} (\f'(a_i)) \bt'_h(\f'(a)) =0$ and applying $\bt'_{h\m}$ we obtain $\sum_i \bt'_{h\m g_i} (\f'(a_i)) \f'(a)
=0.$ Since  $\f(\A')$ is an ideal in $\B'$ the element $\bt'_{h\m g_i} (\f'(a_i)) \f'(a)$ is contained in the algebra  $\bt'_{h\m
g_i}(\f'(\A)) \cap \f'(\A) = \f'(\D_{h\m g_i}) = \f'(\A 1_{h\m g_i}) = \f'(\A) \f'(1_{h\m g_i}),$ whose unity element is
$\f'(1_{h\m g_i}).$ Therefore, using (ii'),

\begin{align*}
& \bt'_{h\m g_i} (\f'(a_i)) \; \f'(a) = \bt'_{h\m g_i} (\f'(a_i)) \; \f'(1_{h\m g_i}) \; \f'(a) =\\&=
\bt'_{h\m g_i} (\f'(a_i)) \; \f' \circ \af_{h\m g_i}(1_{g_i\m h})  \;\f'(a) = \bt'_{h\m g_i} ( \f'(a_i))   \; \bt'_{h\m
g_i} \circ \f' (1_{g_i\m h}) \;  \f'(a) = \\&= \bt'_{h\m g_i} \circ \f'(a_i  1_{g_i\m h}) \; \f'(a) = \f' \circ \af_{h\m g_i}
 (a_i  1_{g_i\m h}) \; \f'(a) =  \f' ( \af_{h\m g_i}  (a_i  1_{g_i\m h}) \; a).
\end{align*}

In a similar fashion we see that $ \bt_{h\m g_i} (\f(a_i))\; \f(a) = \f ( \af_{h\m g_i}  (a_i  1_{g_i\m h}) \; a).$ 
Thus,

\begin{align*}
&     0 = \sum_{i=1}^s \bt'_{h\m g_i} (\f'(a_i)) \f'(a)   =   \sum_{i=1}^s  \f' ( \af_{h\m g_i}  (a_i  1_{g_i\m h}) \; a),
\end{align*}

\noindent which implies    $ \sum_{i=1}^s   \af_{h\m g_i}  (a_i  1_{g_i\m h} \; )a =0. $ Hence $
 0 = \sum_{i=1}^s   \f ( \af_{h\m g_i}  (a_i  1_{g_i\m h}) \; a) =  \sum_{i=1}^s   \bt_{h\m g_i} (\f(a_i))\; \f(a),$ and
applying $\bt_h$ we obtain

\begin{align*}
&      \sum_{i=1}^s \bt_{ g_i} (\f(a_i))\; \bt_h(\f(a)) = 0,
\end{align*}

\noindent for all $a \in \A$. Therefore, the element $\sum_{i=1}^s \bt_{g_i} (\f(a_i))$ annihilates each
$\bt_h(\f(\A)).$ Let $\B_1$ be the  algebra generated by $\cup_{i=1}^s \bt_{g_i} (\f(\A)).$ Then $\sum_{i=1}^s
\bt_{g_i} (\f(a_i)) \in
\B_1$ and by Lemma~\ref{unity} $\B_1$ posses the unity element, say $1_{\B_1}.$ Then     $$\sum_{i=1}^s \bt_{ g_i} (\f(a_i)) =
\sum_{i=1}^s
\bt_{ g_i} (\f(a_i)) \cdot 1_{\B_1} = 0,$$   so that $\sum_{i=1}^s \bt'_{ g_i} (\f'(a_i)) = 0$ implies $\sum_{i=1}^s \bt_{ g_i}
(\f(a_i)) = 0$ and $\phi: \B' \to \B$ is a well-defined homomorphism of algebras.

By symmetry,  $\bt_g(\f(a)) \mapsto \bt'_g(\f'(a)), \;g\in G, a\in \A$ also determine a
well-defined   map $\phi': \B \to \B'.$ Obviously, $\phi' \circ \phi = \phi \circ \phi' = 1$ and, consequently, $\phi$ is an
isomorphism of algebras. It is easily seen that for all $g \in G$ one has $\bt_g \circ \phi = \phi \circ \bt'_g$ and this
yields that  $\bt'$ is   equivalent to $\bt$.

\fim

\end{section}

\begin{section}{Morita equivalence}\label{sec:morita}
Recall from \cite[Section 4.1]{Rowen} that a {\it Morita context}\/ is a
six-tuple $(R,R',M,M',\tau,\tau')$, where:

\medskip
  (a) $R$ and $R'$ are rings,

  (b) $M$ is an $R$-$R'$-bimodule,

  (c) $M'$ is an $R'$-$R$-bimodule,

  (d) $\tau : M\otimes_{R'} M' \to R$ is a bimodule map,

  (e) $\tau' : M'\otimes_{R} M \to R'$ is a bimodule map,

\medskip\noindent  such that
  $$
  \tau(x\otimes x')\;y = x\; \tau'(x'\otimes y),\qquad \forall x,y\in
M,\; x'\in M',$$
  and
  $$\tau'(x'\otimes x)\;y' = x'\; \tau(x\otimes y'),\qquad \forall
x',y'\in M',\; x\in M.$$

By Morita's fundamental results \cite[Theorems 4.1.4 and
4.1.17]{Rowen}, given a Morita context with $\tau$ and $\tau'$ onto,
the categories  of $R$-modules and of $R'$-modules are equivalent.
In this case $R$ and $R'$ are said to be {\it Morita
equivalent}.

Let $\af$ be a partial action of a group $G$ on a unital algebra $\A$ and
suppose that $(\bt,\B)$ is an enveloping action for $(\af,\A)$ such that   $\B$
also has unity. It is
our goal in this section to exhibit an explicit Morita context for the
rings $R=\A\stimes_\af G$ and $R'=\B\stimes_\bt G$, with $\tau$ and
$\tau'$ onto, hence proving these to be Morita equivalent.  This
should be seen as a purely algebraic counterpart to \cite[Theorem
4.18]{AbadieTwo}.

Consider the
linear subspaces $M,N\subseteq \B\stimes_\bt G$ given by
  $$
  M = \left\{\sum_{g\in G}c_g\de_g : c_g\in \A\text { for all }g\right\}
  $$
and
  $$
  N = \left\{\sum_{g\in G}c_g\de_g : c_g\in \bt_g(\A)\text { for all }g\right\}.
  $$

  \begin{prop}
  $M$ is a right ideal and $N$ is a left ideal of
$\B\stimes_\bt G$.
  \end{prop}

\p Let $c\de_g$ be in $M$, where $g\in G$ and $c\in \A$.  If $h\in G$
and $b\in\B$ we have that
  $$
  c\de_g \cdot b\de_h =
  c\bt_g(b)\de_{gh} \in M,
  $$
  because $c\bt_g(b)\in\A$ since $\A$ is an ideal in $\B$.  So
$M$ is a right ideal in $\B\stimes_\bt G$.

Next let $c\de_g$ be in $N$, where $g\in G$ and $c=\bt_g(c')$ with
$c'\in \A$.  If $h\in G$ and $b\in\B$ we have
  $$
  b\de_h \cdot  c\de_g  =
  b\bt_h(c)\de_{hg}  =
  b\bt_{hg}(c')\de_{hg} \in N,
  $$
  because  $\bt_{hg}(\A)$ is an ideal in $\B$.  So  $N$ is a left
ideal.
  \fim\\

In our next result we will view  $\A\stimes_\af G$ as a
subalgebra of $\B\stimes_\bt G$ as in  (\ref{envelopecross}).

  \begin{prop}
  $(\A\stimes_\af G)M \subseteq M$ and $N(\A\stimes_\af
G)\subseteq N$ so that $M$ may be viewed as a left $\A\stimes_\af
G$-module and $N$ may be viewed as a right $\A\stimes_\af G$-module.
  \end{prop}

\p Let $c\de_g$ be in $M$, where $g\in G$ and $c\in \A$, and let
$a\de_h\in \A\stimes_\af G$, where $h\in G$ and $a\in\D_h$.  Then
  $$
  a\de_h \cdot c\de_g =
  a \bt_h(c)\de_{hg} \in M,
  $$
  because $\A$ is an ideal in $M$.    So $M$ is a left $\A\stimes_\af G$-module.

Next let $c\de_g$ be in $N$, where $g\in G$ and $c=\bt_g(c')$ with
$c'\in \A$, and let $a\de_h\in \A\stimes_\af G$, where $h\in G$ and $a\in\D_h$.  Then
  $$
  c\de_g \cdot a\de_h  =
  c\bt_g(a)\de_{gh} =
  \bt_g(c'a)\de_{gh} = \cdots
  $$
  Since $\D_h$ is an ideal in $\A$ we have that $c'a\in\D_h$ and hence
$c'a=\af_h(x)$ for some $x\in \D_{h\m}$.  So the above equals
  $$
  \cdots =
  \bt_g(\af_h(x))\de_{gh} =
  \bt_g(\bt_h(x))\de_{gh} =
  \bt_{gh}(x)\de_{gh} \in N.
  $$
  So $N$ is a right $\A\stimes_\af G$-module.
  \fim\\

So far we have therefore seen that $M$ is an ($\A\stimes_\af
G$)-($\B\stimes_\bt G$)-bimodule, whereas $N$ is a ($\B\stimes_\bt
G$)-($\A\stimes_\af G$)-bimodule.  These give the third and fourth
components of the Morita context we are looking for.  Next we need to
find $\tau$ and $\tau'$.

Given linear subspaces $X$ and $Y$ of an algebra we will denote by
$XY$ the linear span of the set of products $xy$ with $x\in X$ and $y\in Y$.  This notation will be used in our next:

  \begin{prop}\label{epic}
  $MN=\A\stimes_\af G$ and
$NM=\B\stimes_\bt G$.
  \end{prop}

  \p
  Let $c_1\de_g\in M$ and $c_2\de_h\in N$, where  $g,h\in G$,
$c_1\in \A$, and $c_2=\bt_h(c_2')$ with  $c_2'\in\A$. Then
  $$
  c_1\de_g \cdot c_2\de_h =
  c_1\bt_g (c_2)\de_{gh} =
  c_1\bt_{gh} (c'_2)\de_{gh}
  \in \A\stimes_\af G,
  $$
  because  $c_1\bt_{gh} (c'_2)\in \A\cap\bt_{gh}(\A) = \D_{gh}$.
  This shows that $MN\subseteq \A\stimes_\af G$.

  Given $c\in\D_h$ observe that
$c\in\af_h(\D_{h\m})\subseteq\bt_h(\A)$ so that $c\de_h\in N$.
Letting $1_\A$ be the unit of $\A$ we have that $1_\A\de_e\in M$ and
  $$
  1_\A\de_e\cdot c\de_h =
  1_\A c \de_h =
  c \de_h.
  $$
  So $c\de_h\in MN$.
  Since $c$  is arbitrary we conclude that $\A\stimes_\af G\subseteq MN$.

  Let  $g,h\in G$ and let $c\in \A$.  Then $\bt_g(c)\de_g\in N$
and $1_\A\de_{g\m h}\in M$.  Moreover
  $$
  \bt_g(c)\de_g\cdot 1_\A\de_{g\m h} =
  \bt_g(c 1_\A)\de_{h} =
  \bt_g(c)\de_{h}.
  $$
  Since $\cup_{g \in G} \bt_g(\A)$ generates $\B$ we conclude that
$\B\de_h\subseteq NM$ and hence $NM=\B\stimes_\bt G$.
  \fim

The above result suggests that we define
  $$
  \tau : m\otimes n \in M \otimes_{_{\B\stimes_\bt G}} N \longmapsto
mn \in \A\stimes_\af G
  $$
and
  $$
  \tau' : n\otimes m \in N \otimes_{_{\A\stimes_\af G}} M \longmapsto
nm \in \B\stimes_\bt G.
  $$
  It is then easy to see that the six-tuple
  $(\A\stimes_\af G, \B\stimes_\bt G, M, N, \tau,\tau')$ is a Morita
context.  This leads us to the main result of this section:

\begin{teo} Let $\af$ be a partial action of a group $G$ on a unital algebra
$\A$ and suppose that $(\bt,\B)$ is an enveloping action for
$(\af,\A)$ such that $\B$ is also unital.  Then $\A\stimes_\af G$ and $\B\stimes_\bt G$ are Morita
equivalent.  \end{teo}

\p By (\ref{epic}) we have that both $\tau$ and $\tau'$ are onto.  The
conclusion then follows from \cite[Theorem 4.1.17]{Rowen}.  \fim

\end{section}

\begin{section}{Partial representations and partial actions}\label{sec:reprcross}

        In this section we use crossed products to relate partial actions with  {\it  partial representations} of groups.

\begin{defi} A {\it partial representation} of a group
$G$ into a unital $K$-algebra $\B$ is a map  \[ \pi : G \to \B, \]
  which sends the unit  element of the group to the unity element of $\B$, such that
for all $g,h\in G$ we have

$$\pi(g)\pi(h)\pi(h^{-1})=\pi(gh)\pi(h^{-1})\ \ \text{ and }\ \ \pi(g^{-1})\pi(g)\pi(h)=\pi(g^{-1})\pi(gh).\\$$

\end{defi}

In particular, if $\B$ is the algebra of linear transformations ${\rm End}(V)$ of a vector space $V$ over a field
$K$ then we have a partial representation of $G$ on the vector space $V$. \\

\begin{lema}\label{map} Let $\af = \{\af_g : \D_{g\m} \to \D_g \; (g \in G)\}$ be a partial action of $G$ on an algebra
$\A$ such that each $\D_g \; (g\in G)$ is a unital algebra with unity $1_g$. Then the map
$\pi_{\af}: G \ni g \mapsto 1_g \de_g \in \A \stimes_{\af} G$ is a partial representation.
\end{lema}

\p Obviously,  $1 \mapsto  1 \cdot \de_1$, the identity element of  $\A \stimes_{\af} G.$  We see that

\begin{align*}
&(1_{g\m} \de_{g\m} \cdot 1_g \de_g) \cdot 1_h \de_h =\af_{g\m}(\af_g(1_{g\m})1_g)  \de_1 \cdot 1_h \de_h =
\af_{g\m} (1_g^2) \de_1 \cdot 1_h \de_h = \\ & =
\af_{g\m} (1_g) \de_1  \cdot 1_h \de_h =  1_{g\m} \de_1 \cdot 1_h \de_h = 1_{g\m} 1_h \de_h.
\end{align*}

\noindent On the other hand, using (\ref{triviality}), we have

\begin{align*}
& 1_{g\m} \de_{g\m} \cdot 1_{gh} \de_{gh}  = \af_{g\m}(\af_g(1_{g\m}) 1_{gh})  \de_h = \af_{g\m}(1_g 1_{gh}) \de_h =
   1_{g\m} 1_h \de_h.
\end{align*}

\noindent Thus $1_{g\m} \de_{g\m} \cdot 1_g \de_g \cdot 1_h \de_h =  1_{g\m} \de_{g\m} \cdot 1_{gh} \de_{gh}$ for all $g,h \in
G$ and one similarly verifies that  $1_{g}  \de_g \cdot 1_h \de_h  \cdot 1_{h\m} \de_{h\m} = 1_g 1_{gh} \de_g =  1_{gh}
\de_{gh} \cdot 1_{h\m} \de_{h\m}$ $(g,h \in G),$ which shows that $\pi_{\af}: G \to \A \stimes_{\af} G$ is a partial representation.
\fim \\

\begin{defi} Two partial representations $ \pi : G \to \B$ and $ \pi' : G \to \B'$ are equivalent if there is an
isomorphism $\f:\B' \to \B$ such that $$ \pi(g) = \f (\pi' (g)) $$ for all $g \in G.$
\end{defi}

\begin{remark}\label{remark1} It is easily verified that the map $\af \mapsto \pi_{\af}$ sends equivalent partial actions into
equivalent partial representations.
\end{remark}

Let $ \pi : G \to \B $ be a partial representation of a group $G$ into a unital $K$-algebra $\B$. Then
by~(2),~(3) of~\cite{DEP}
the elements $\e_g=\pi(g)\pi(g\m)$ $(g\in G)$ are commuting
idempotents such that

\begin{equation}\label{1star}
\pi(g)\e_h=\e_{gh}\pi(g),\quad \e_h\pi(g)=\pi(g)
\e_{g\m h}.
\end{equation}

Let $\A$ be the subalgebra of $\B$ generated by all the $\e_g$ $(g\in G)$ and for a fixed $g \in G$ set  $\D_g = \e_g \A$.

\begin{lema}\label{action} The maps $\af^{\pi}_g  :{\D}_{g\m} \to  {\D}_g$ $(g \in G),$
defined by  $\af^{\pi}_g (a) = \pi(g) a \pi(g\m)$ $(a \in {\D}_{g\m}),$  are  isomorphisms of $K$-algebras, which determine a
partial action $ \af^{\pi}  $ of $G$ on $\A$.
\end{lema}

\p Write for simplicity $\af^{\pi} = \af.$ We observe first that $\A$ is invariant with respect to the map $a \mapsto \pi(g) a
\pi(g\m).$ Obviously,
$\A$ is spanned  by elements $a = \e_{h_1} \ldots \e_{h_s}$ with $h_1, \ldots, h_s \in G.$  We see that

%\begin{align*}

%&

$$  \pi(g)  \e_{h_1} \ldots \e_{h_s} \pi(g\m) =    \pi(g)  \pi(g\m) \e_{gh_1} \ldots \e_{gh_s} = \e_g \e_{gh_1}
\ldots \e_{gh_s} \in \A.$$

%= \\ & =

%&
%\end{align*}

\noindent Thus $\pi(g) \A \pi(g\m) \subseteq \A.$ It is easily verified that this map sends $\e_{g\m}$ to $\e_{g}.$
Since  $\D_{g}$ is generated in $\A$ by the idempotent $\e_g$, we obtain a $K$-map $\af_g  :{\D}_{g\m} \to  {\D}_g.$
Moreover, it is a homomorphism of algebras. Indeed, taking $a, b \in \D_{g\m}$ we see that

\begin{align*}
& \af_g (a) \af_g(b) = \pi(g) a \pi(g\m) \pi(g) b \pi(g\m) = \pi(g) a \e_{g\m} b \pi(g\m) = \pi(g) ab \pi(g\m) = \af_g(ab),
\end{align*}

\noindent as $\e_{g\m}$ is the identity element of $\D_{g\m}.$ Then  $\af_{g\m}  :{\D}_g \to  {\D}_{g\m}$ is also a
homomorphism of algebras and it is
easily seen that $\af_g \circ \af_{g\m}$ and  $\af_{g\m} \circ \af_g$ are identity maps.

Taking $g, h \in G$ write an element  $a \in {\D}_h \cap {\D}_{g\m}$
as $ a = \e_h \e_{g\m} b$ with $b \in \A$. Then

$${\af}\m_h(a) = \pi(h\m) \e_h \e_{g\m} b \pi(h) = \pi(h\m) \e_h \pi(h) \e_{h\m g\m} b' =  \e_{h\m}  \e_{ h\m g\m } b' \in
\D_{(gh)\m},$$

\noindent with $b' \in \A.$ Hence  $\D_{(gh)\m} \supseteq {\af}\m_h({\D}_h \cap {\D}_{g\m}),$ which is condition $(ii)$ of
Definition~\ref{def1}.

Finally, for $a \in  {\af}\m_h ({\D}_h \cap {\D}_{g\m})$, we have $ \e_{h\m} a= a \e_{h\m} = a,$ as $a \in {\D}_{h\m}$ and,
thus,

\begin{align*}
& \af_g \circ \af_h (a) = \pi(g) \pi(h) a \pi(h\m) \pi(g\m) = \pi(g) \pi(h) a  \e_{h\m} \pi(h\m) \pi(g\m) = \\
& = \pi(g) \pi(h) a  \e_{h\m} \pi(h\m g\m) = \pi(g) \pi(h)   \e_{h\m} a \pi(h\m g\m) = \\
& = \pi(g h)   \e_{h\m} a \pi(h\m g\m) = \af_{gh} (a),
\end{align*}

\noindent which gives $(iii)$ of Definition~\ref{def1}. Since $(i)$ is obvious, the Lemma is proved. \fim\\

\begin{remark}\label{remark2} As in the previous case it is readily   seen that  $\pi \mapsto \af^{\pi}$ also preserves
equivalence relations.
\end{remark}

\begin{prop}\label{hom1} Let $\af = \{\af_g : \D_{g\m} \to \D_g \; (g \in G)\}$ be a partial action of $G$ on an algebra
$\A$ such that each $\D_g \; (g\in G)$ is a unital algebra with unity $1_g$. Let further $\A'$ be the subalgebra of
$\A\stimes_{\af}G$ generated by all $1_g  \de_1$ $(g \in G).$ Then the map $\f_{\af}: \A' \ni 1_g \de_1 \mapsto 1_g \in \A$ is a
monomorphism such that $\f_{\af} \circ {\af}^{\pi_{\af}}_g = \af_g \circ \f_{\af}$ for each $g \in G.$ In particular, if $\A$
is generated by the elements $1_g$ $(g \in G),$ then the partial actions $\af^{\pi_{\af}}$ and $\af$ are equivalent.
\end{prop}

\p Clearly, $\f_{\af}$ is the restriction of the isomorphism $\A \cdot \de_1 \ni a \de_1  \mapsto a \in \A $ and thus is
a monomorphism of $\A'$ into $\A$.

According to Lemma~\ref{map} we have the partial representation $\pi_{\af}: G \ni g \mapsto 1_g \de_g \in \A \stimes_{\af} G$
which  by Lemma~\ref{action} induces a partial action on the subalgebra of $\A \stimes_{\af} G,$ generated by the elements
$\pi_{\af}(g) \pi_{\af}(g\m).$ Because $$\pi_{\af}(g) \pi_{\af}(g\m) = 1_g \de_g 1_{g\m} \de_{g\m} = \af_g(\af_{g\m}(1_g)
1_{g\m}) \de_1  = 1_g \de_1,$$ this subalgebra is exactly $\A'.$   An arbitrary element of $\A'$ can be written as $a'
\de_1$ where $a'$ belongs to the image $Im \f_{\af}$ of $\f_{\af}.$ The partial action
$\af^{\pi_{\af}}$ is given by the isomorphisms $\af^{\pi_{\af}}: \D'_{g\m} \ni a' \de_1 \mapsto 1_g \de_g \cdot a' \cdot
1_{g\m}
\de_{g\m} \in
\D'_{g},$ where
$\D'_{g}$ is the ideal of $\A'$ generated by $1_g \de_1.$ For $a' \de_1 \in \D'_{g\m}$ we see that

\begin{align*} & \f_{\af} ( {\af}^{\pi_{\af}}_g (a' \de_1 )) = \f_{\af} ( 1_g \de_g a'   1_{g\m} \de_{g\m}) =
 \f_{\af} (\af_g( \af_{g\m}(1_g) a' 1_{g\m})\de_1) = \\& \f_{\af} (\af_g( a' 1_{g\m})\de_1) =  \af_g( a') =
\af_g(\f_{\af}(a' \de_1)),
\end{align*}
\noindent so that $\f_{\af} \circ {\af}^{\pi_{\af}}_g =  \af_g \circ \f_{\af}$ for all $g \in G.$ If $\A$
is generated by the elements $1_g$ $(g \in G),$ then clearly $\f_{\af} : \A' \to \A$ is an isomorphism which gives the
equivalence of the partial actions $\af^{\pi_{\af}}$ and $\af.$ \fim

\begin{prop}\label{hom2} Let $\pi:G \to \B$ be a partial representation and suppose that the subalgebra $\A \subseteq \B$ and
the partial action $\af^{\pi}$ of $G$ on $\A$ are  as in Lemma~\ref{action}. Then the map
$\phi_{\pi} : \A \stimes_{\af^{\pi}} G \to  \B,$  defined by  $\phi_{\pi}( \sum_{g \in G} a_g  \de_g) =  \sum_{g \in G} a_g  \pi(g),
$   is  a homomorphism of $K$-algebras such that $\phi_{\pi} \circ \pi_{\af^{\pi}} = \pi.$ In particular, if $\phi_{\pi}$ is an
isomorphism, then  the partial representations $\pi$ and $\pi_{\af^{\pi}}$ are equivalent.
\end{prop}

\p Obviously, the elements $a_g  \de_g$ $(g \in G, a_g \in \D_g)$ span $\A \stimes_{\af^{\pi}} G$. Since $\e_g$ is the identity
element of $\D_g$ and  $\e_g \pi(g h) = \e_g \pi(g) \pi(h),$ we have

\begin{align*}
& \phi_{\pi}( a_g \de_g  \cdot b_h \de_h) = \\ & =  \phi_{\pi}( \af^{\pi}_g ( \af^{\pi}_{g\m}(a_g) b_h)  \de_{gh} ) =
\af^{\pi}_g ( \af^{\pi}_{g\m}(a_g) b_h)
\pi(gh)  =   \af^{\pi}_g ( \af^{\pi}_{g\m}(a_g) b_h) \e_g \pi(gh)  = \\&=  \af^{\pi}_g ( \af^{\pi}_{g\m}(a_g) b_h) \e_g \pi(g)
\pi(h)  = \af^{\pi}_g ( \af^{\pi}_{g\m}(a_g) b_h)  \pi(g) \pi(h) = \\&=
 \pi(g) ( \af^{\pi}_{g\m}(a_g) b_h) \pi(g\m)  \pi(g) \pi(h) =   \pi(g)  \af^{\pi}_{g\m}(a_g) b_h \e(g\m)   \pi(h) = \\&=
\pi(g)  \af^{\pi}_{g\m}(a_g) b_h \pi(h),
\end{align*}

\noindent as $\af^{\pi}_{g\m}(a_g) b_h \in \D_{g\m}.$  Continuing this calculation, we see that

\begin{align*}
& \pi(g)  \af^{\pi}_{g\m}(a_g) b_h \pi(h) =   \pi(g) \pi(g\m) a_g \pi(g) b_h \pi(h) =   \e(g) a_g  \pi(g) b_h \pi(h)  = \\&=
a_g  \pi(g) b_h \pi(h) =
\phi_{\pi}(a_g \de_g)  \cdot \phi_{\pi}( b_h \de_h),
\end{align*}

\noindent using once more that $\e_g$ is the identity of $\D_g$. Thus, $\phi_{\pi} : \A \stimes_{\af^{\pi}} G \to  \B$ is a
homomorphism of algebras. It is easily seen that $\phi_{\pi} \circ \pi_{\af^{\pi}} (g) = \phi_{\pi}( \e_g \de_g) = \e_g \pi(g)
= \pi(g) \pi(g\m) \pi(g) = \pi(g),$  for each $g \in G,$ which shows that $\phi_{\pi} \circ \pi_{\af^{\pi}} = \pi.$ If
$\phi_{\pi}$ is an isomorphism, then it  obviously follows that  $\pi$ and $\pi_{\af^{\pi}}$ are equivalent. \fim\\

Now we can show that the partial group algebras can be naturally
endowed with a crossed product structure. We remind the reader that
the partial group algebra $\kpar(G)$ of a group $G$ over $K$ is the
semigroup algebra $K S(G),$ where $S(G)$ is the universal semigroup,
generated by the symbols $\{ [g] : g \in G \}$ subject to relations:\\

$ a) [g\m][g][h] = [g\m][gh];$

$ b) [g][h][h\m]=[gh][h\m] \; \; \; (g,h \in G);$

$ c) [1]=1,$\\

\noindent where  $1$ also denotes   the identity
element of $S(G)$. It is easily seen that $\tilde\pi: G \ni g \mapsto [g] \in \kpar(G)$ is a partial representation and
$\kpar(G)$ has the following universal property: for every partial representation $\pi$ of $G$ into a unital algebra $\B$
there exists a unique homomorphism of algebras $\psi: \kpar(G) \to \B$ such that $\pi = \psi \circ \tilde\pi$, and, conversely,
for each algebra homomorphism $\psi: \kpar(G) \to \B$ the map $\psi \circ \tilde\pi$ is a partial representation of $G$ into
$\B.$

\begin{teo}\label{pgr} Let $\A$ be the subalgebra of $\kpar(G)$ generated by the elements $\tilde\e_g = \tilde\pi(g)
\tilde\pi(g\m) = [g] [g\m]$.  Then 
the homomorphism $\phi_{\tilde\pi}: \A
\stimes_{\af^{\tilde\pi}} G \to \kpar(G)$ associated to the partial
representation $\tilde\pi: G \ni g \mapsto [g] \in \kpar(G) $ is
actually an isomorphism.\end{teo}

\p By Proposition~\ref{hom2}, $\phi_{\tilde\pi} : \A \stimes_{\af^{\tilde\pi}} G \to \kpar(G) $ is a homomorphism and we shall show
that it has an inverse.  Recall that $\D_g$ is spanned by the elements of form $\tilde\e_g
\cdot \tilde\e_{h_1}  \cdot \ldots \cdot  \tilde\e_{h_s}.$ By the universal property of partial group
rings,  the partial representation
$\pi_{\af^{\tilde\pi}}: G \ni g \mapsto \tilde\e_g \de_g \in \A \stimes_{\af^{\tilde\pi}} G$
gives a homomorphism $\psi: \kpar(G) \to \A \stimes_{\af^{\tilde\pi}} G$ such that   $\psi( [g])  = \tilde\e_g \de_g.$

Obviously, $\phi_{\tilde\pi} \circ \psi [g] = \phi_{\tilde\pi} (\tilde\e_g \de_g) = \tilde\e_g [g] = [g][g\m][g] = [g]$
for each $g \in G,$ so that $\phi_{\tilde\pi} \circ \psi$ is the identity map. On the other hand, since $\tilde\e_g
\de_g \tilde\e_{g\m} \de_{g\m} = \tilde\e_g \de_1,$ we have

\begin{align*} & \psi \circ \phi_{\tilde\pi} ((\tilde\e_g \cdot \tilde\e_{h_1}  \cdot \ldots \cdot  \tilde\e_{h_s})\de_g)=
\psi  ((\tilde\e_g \cdot \tilde\e_{h_1}  \cdot \ldots \cdot  \tilde\e_{h_s})[g])=\\&=
\psi([g][g\m][h_1][h\m_1]\ldots[h_s][h\m_s][g]) =\\ & =  (\tilde\e_g \de_g \tilde\e_{g\m} \de_{g\m}) (\tilde\e_{h_1} \de_{h_1}
\tilde\e_{h\m_1} \de_{h\m_1}) \ldots (\tilde\e_{h_s} \de_{h_s} \tilde\e_{h\m_s} \de_{h\m_s}) \tilde\e_g \de_g =\\&=
\tilde\e_g \de_1 \tilde\e_{h_1} \de_1 \ldots \tilde\e_{h_s} \de_1 \tilde\e_g \de_g = \tilde\e_g^2  \tilde\e_{h_1}
\ldots \tilde\e_{h_s}  \de_g = \tilde\e_g  \tilde\e_{h_1}  \ldots \tilde\e_{h_s}  \de_g.
\end{align*}

\noindent Thus $\psi \circ \phi_{\tilde\pi}$ is also the identity map and consequently $\psi$ is the inverse of
$\phi_{\tilde\pi}.$ \fim \\

The next fact helps to identify algebras as crossed products.

\begin{prop}\label{monomorphism} Let   $\af = \{\af_g : \D_{g\m} \to \D_g \; (g \in G)\}$ be a partial action of $G$ on an
algebra $\A$ such that for each $g \in G$ the right annihilator of  $\D_g$ in  $\D_g$ is zero. Suppose that
$\f : \A \stimes_{\af} G \to \B$ is a homomorphism of algebras whose restriction to $\A$  is injective. If there exists
a
$K$-linear transformation $E: \B \to \B$ such that for every $g \in G, \; a \in \D_g$

\begin{align*}
 E(\f(a \de_g)) =
\begin{cases} \f(a \de_g) & \text{if   } g=1,\\
               0          & \text{if   } g \neq 1.
\end{cases}
\end{align*}

\noindent then $\f$ is injective.
\end{prop}

\p Suppose that $\sum_{g \in G} a_g \de_g$ belongs to the kernel $Ker \f$ of $\f$. Then for each $b \in \D_{h\m}$ we have
$\sum a_g
\de_g
\cdot b
\de_{h\m}
\in Ker \f$ and

\begin{align*} 0 = \f ( \sum_{g \in G} a_g \de_g \cdot  b \de_{h\m} ) = \f ( \sum_{g \in G} \af_g ( \af_{g\m} (a_g) b) \de_{g
h\m} ).
\end{align*}

\noindent Applying $E$ we eliminate all terms with $g \neq h$ so that $$ \f (  \af_h ( \af_{h\m} (a_h) b) \de_1 ) =0.$$
Because $\f$ is injective on $\A$ $(= \A \de_1)$ we obtain that $$0 = \af_h(\af_{h\m}(a_h) b) = a_h \af_h(b)$$ for every $b
\in \D_{h\m}.$  Since $\af_h(b)$ runs over $\D_h$ and the annihilator of  $\D_h$ in $\D_h$ is zero, we conclude that $a_h = 0
$ for each $h
\in G.$ This yields that $\f$ is a monomorphism. \fim\\

We apply the above results to matrix partial representations. First we recall some facts obtained in \cite{DEP} and
\cite{DZh0} (see also \cite{DZh}). Let $A$ be a subset of a group $G$ which contains $1$ and let $H$ be the stabilizer of $A$
in
$G$ with respect to multiplication from left, that is
$H = St(A) = \{h \in G: hA=A\}.$ Then
$A$ is a union of right cosets of $H.$ Suppose that $A$ is a finite union of such cosets, i.e. the index $(A:H)$ of $H$ in
$A$ is finite, say $n$. Then left multiplication of
$A$ by the elements
$g$ with $g\m \in A$ give a finite number of distinct sets $A = A_1,
A_2, \ldots, A_n,$ each of which contains $1.$
  Consider the category $\G$ whose objects are the sets $A_1, \ldots,
A_n$ and the morphisms are as follows.  If $A_i$ and $A_j$ are
given and there is some $g\in G$ such that  $gA_i = A_j$ (in which
case $g\m$ necessarily belongs to $A_i$) then Hom$(A_i,A_j)$ is the
set of all pairs $(A_i,g)$ with $g$ as above.  Otherwise
Hom$(A_i,A_j)$ is the empty set.
  A
morphism can be viewed as  left multiplication of $A_i$ by $g,$ where $g\m \in A_i.$  The product $(A_i,g) \circ (A_j,g')$ is
defined in
$\G$ if and only if $g'A_j = A_i$ and if so, is equal to $(A_j, gg').$ The identity morphism of
$A_i$ is clearly
$(A_i,1)$ and, because each morphism has an inverse, $\G$ is a groupoid. The subgroup $H$ is the isotropy group of $A$ and
the isotropy group of $A_i$ is conjugate to $H$ in $G.$ Identify $\G$ with the set of its morphisms. Then the groupoid algebra
$K \G$ is the $K$-vector  space with basis $\G$ endowed with the following multiplication:

\[\gm_1\cdot\gm_2=\left\{
\begin{array}{ll}\gm_1\gm_2,&\text{if the composite morphism}\ \gm_1\gm_2
 \text{ exists in }\G,
\\ 0,  &  \text{otherwise}.
\end{array}
\right.\]

By \cite[Proposition 3.1]{DEP}, the groupoid algebra $K \G$ is isomorphic to the algebra $M_n(K H)$ of $n\times n$-matrices
over the group algebra $K H$ and an isomorphism $\tau: K \G \to M_n(K H)$ can be obtained as follows. For each $A_i$ fix an
element $g_i \in G$  such that $g_i A = A_i.$ If $g A_i = A_j$, we see that $g\m_j g g_i = h$ for some $h \in
H.$ Then the map $\tau$ is given by $(A_i,g)  \mapsto e_{j,i}(h),$ where $e_{j,i}(h)$ is the elementary
matrix whose unique non-zero entry is $h\in H$, which is placed in the  intersection of the $j$'th row and $i$'th column.
Consider also the map $\lb:G \to K\G$, defined by

\begin{equation*}\label{eq:lpar}
\lb(g)=\left\{
\begin{array}{ll} \sum_{\substack{ A_i\ni g^{-1}}}(A_i,g), &\text{if}\ g\m \in \cup_i A_i,\\
0, & \text{otherwise}.
\end{array}
\right.
\end{equation*}
According to  \cite{DZh} $\lb$ is a partial representation of $G$ into $K \G.$ If  $\s': M_n(K H) \to M_{nm}(K)$ is the
representation extended from an irreducible representation $\s: H \to M_m(K)$ then $ \s' \circ  \tau \circ \lb : G \to
M_{nm}(K)$ is an irreducible partial representation and, moreover, every irreducible  partial matrix $K$-representation of $G$
is equivalent to a partial representation of the form  $ \s' \circ  \tau \circ \lb$ for some subset $A
\ni 1$ of $G$ and some irreducible matrix $K$-representation $\s$ of the isotropy group $H$ of $A.$ Thus the ``partial part''
of an irreducible finite dimensional partial $K$-representation is given by the partial representation $\tau \circ \lb : G
\to  M_n(K H)$. The partial representations of the form $\tau \circ \lb$ obtained from the subsets $A \ni 1$
shall be called the {\it elementary partial representations} of $G.$  Different choices of the elements $g_i$ with $A_i = g_i
A$ give rise to equivalent partial representations. By the {\it equivalence class} of an elementary partial representation we
mean the set of all partial representations of $G$ which are equivalent to the given elementary partial representation. An
elementary partial representation has {\it trivial isotropy}   if the corresponding isotropy group
$H$ is $1.$

  We shall see that an elementary partial representation $\pi: G \to
M_n(KH)$ gives rise to a certain partial action of $G$
on the diagonal subalgebra $diag(K,\ldots,K)$ of the full matrix algebra  $M_n(K).$ This diagonal algebra is obviously
isomorphic to
$K^n,$  the $n$'th direct power of $K,$ so we shall speak about partial actions on $K^n.$ Given a partial action
$\af$ of $G$ on $diag(K,\ldots,K) \cong K^n$ we see that  each $\D_g$ is an algebra with unity
$1_g$ which is a sum of some minimal idempotents $e_{i,i}(1).$ Set $A_{i}(\af) = \{ g \in G : (1_g)_{i,i} \neq 0 \},$ where
$(1_g)_{i,i}$ denotes the $i$'th diagonal entry of $1_g.$

\begin{teo}\label{correspondence1}  For a fixed  $n>0$ and a fixed subgroup $H$ of a group $G$ the  maps
 $ \pi \mapsto \af^{\pi}$  and $\af \mapsto \pi_{\af}$   establish a one-to-one correspondence between the equivalence
classes of  the elementary partial representations  $\pi: G \to M_n(KH)$   and  the
equivalence classes of the  partial actions $\af$ of  $G$ on $K^n$ with $St(A_{1}(\af)) =H$ and $(A_{1}(\af) : H) = n$.
\end{teo}

\p  Let $\pi: G \to M_n(KH)$ be an elementary partial representation and $\A$ the subalgebra, generated by all $\e_g$
$(g \in G).$ Then there is a subset $A \ni 1$ in $G$ with isotropy group $H$ such that  $H$  has index $n$ in $A$ and
$\pi$ is of form $\tau \circ \lb $ with $\lb$ and
$\tau$ described above. Thus $\tau$ is determined by a fixed choice of the elements $g_i \in G$ with $A_i=g_i A.$ It is easily
seen that each matrix $\pi(g)$ is ``monomial over
$H$'', i.e. each row and each  column  of $\pi(g)$ contains at most one non-zero entry, which is an element of $H$ (observe
that  zero rows and zero columns are allowed). It follows that $\e_g = \sum_{\substack{ A_i\ni g}} e_{ii}(1)$.  For every $i,
k \in
\{1,\ldots ,n\}$ with $ k \neq i$ there exists $g' \in A_i$ such that $g' \notin A_k$ and we see that the $k$'th diagonal
entry of the matrix
$\e_g
\e_{g'}$ is zero. Hence for a fixed $i$ we have $\prod_{g \in A_i} \e_g  = e_{ii}(1)$ and consequently  $\A$ coincides
with the diagonal subalgebra $diag(K,\ldots,K)$ of $M_n(K).$  Thus by  Lemma~\ref{action},  $\af^{\pi}$ is a
partial action of $G$ on $\A \cong K^n$ and we obviously see that $A_1(\af^{\pi}) = A_1 =A$ so that $St(A_1(\af^{\pi})) = H$
(note that $\e_g = 1_g$).

Next we want to show  that  $\phi_{\pi} : \A \stimes_{\af^{\pi}} G \to  M_n(KH),$  defined in
Proposition~\ref{hom2}, is an isomorphism. Observe first that $\phi_{\pi}$ is an epimorphism.   Indeed, we saw already that the
diagonal subalgebra $diag(K,\ldots,K)$ is contained in the image  $Im (\phi_{\pi})$ of $\phi_{\pi}.$ Moreover, for $g \in G$
we have $\D_g = \oplus_{A_i \ni g} K e_{i,i}(1)$. Fix arbitrarily $h \in H$,  $i,j \in \{1,\ldots,n\}$ and take $g = g_i h
g_j\m$. Then $g \in A_i$ and thus
$e_{i,i}(1) \in \D_g$ so that $ e_{i,i}(1) \de_g$ is an element of $\A \stimes_{\af^{\pi}} G$. We see that the unique non-zero
entry in the $i$'th  row of $\pi(g) = \pi(g_i h g_j\m)$ is $h$ and it is placed in the intersection with the $j$'th
column. Hence $\phi_{\pi} ( e_{i,i}(1) \de_g ) =  e_{i,i}(1) \pi(g) = e_{i,j}(h) \in Im (\phi_{\pi}),$ and because this holds
for  all $h \in H$,  $i,j \in \{1,\ldots,n\},$ it follows that   $Im (\phi_{\pi}) = M_n(KH).$

For an element $x = \sum_{h \in H} x(h)h \in K H$ set $trx = x(1)$ and consider the linear transformation $E: M_n(KH) \ni
(x_{i,j}) \mapsto (tr x_{1,1},
\ldots, tr x_{n,n}) \in diag(K,\ldots,K) \subseteq M_n(K).$  It is easily seen that if $g \neq 1$ then the only possible
diagonal entries of $a \pi(g)$ $(a \in \D_{g})$ are non-identity elements of $H$ multiplied by scalars from $K.$  Thus
$E(\phi_{\pi}(a \de_g)) = 0$ for all $1 \neq g \in G, a \in \D_g$ and  moreover $E( \phi_{\pi}(a \de_1)) = \phi_{\pi}
(a \de_1),$ for each $a \in \A.$  Since the restriction of
$\phi_{\pi}$ to $\A$ is injective, it follows from Proposition~\ref{monomorphism} that $\phi_{\pi}$ is a monomorphism. Hence
$\phi_{\pi}$ is an isomorphism and by Proposition~\ref{hom2} the partial representations $\pi$ and $\pi_{\af^{\pi}}$ are
equivalent.

Suppose now that $\af = \{\af_g : \D_{g\m} \to \D_g \; (g \in G)\}$ is a  partial action of $G$ on $K^n$ such that
$St(A_{1}(\af)) =H$ and the index of $H$ in $A_{1}(\af)$ is $n$. We are going to construct an elementary partial representation
whose $\e_g$'s will be exactly the $1_g$'s.  Suppose for a moment that $G$ is finite
and let $S$ be a subset of $G$. One of the important working tools in \cite{DEP} is the product
$\prod_{g\in S}\e_g \prod_{g \notin S}(1-\e_g)$. Substituting $S$ and $G \setminus S$ by certain finite subsets it is
possible to use this kind of product also for infinite $G$'s, as it is done in the proof of Theorem 1 of \cite{DZh}. We
proceed by adopting these ideas to our situation. Let $G$ be again arbitrary and $S \ni 1$ be a subset of $G$. Each $1_g$ is an
idempotent in
$K^n$ and because this algebra has only a finite number of idempotents, we can choose finite subsets $S' \subset S$ and $S''
\subset G \setminus S$ such that  $\{ 1_g : g \in S \} = \{ 1_g : g \in S' \}$ and $\{ 1_g : g \in G \setminus S \} = \{ 1_g :
g \in S'' \}.$ Set $$f_S = \prod_{g\in S'}1_g \prod_{g \in S''}(1-1_g).$$ Write for simplicity $A_i = A_{i}(\af).$ It follows
from the definition of the $A_{i}(\af)$'s that

\begin{equation}\label{star2}
1_g = \sum_{A_i \ni g } e_{i,i}(1).
\end{equation}

We see that $f_{A_i} \neq 0$ for each $i = 1, \ldots, n,$ as $(f_{A_i})_{ii} = e_{i,i}(1).$ Moreover,

\begin{equation}\label{star3}
f_S \neq 0 \mbox{ if and only if } S=A_i \mbox{ for some } 1 \leq i \leq n.\\
\end{equation}

\noindent Indeed, suppose that $S \neq A_j$ for all $j = 1, \ldots, n$ and fix $i \in  \{1,\ldots,n\}$ arbitrarily. Then
either
 there exists   $g \in S$ with $g \notin A_i$ or there is an element $ t \in A_i$ such that $t \notin S.$ In the
first case $(1_g)_{ii} = 0$ and hence $(f_S)_{ii}=0.$   In the second, $(f_S)_{ii}=0$ because $(1_t)_{ii} = 1$ implies
$(1-1_t)_{ii}=0,$ where $1$ is the identity element of $K^n = diag(K,\ldots,K)$. Since $i$ is arbitrary, it follows that
$f_S=0.$

Next we compute the effect of $\af_g$ on $f_S$ where $g\m \in S.$ It follows from the definition of $f_S$ that taking larger
$S'$ and $S''$ does not alter $f_S.$ Thus we can  choose $S'$ and $S''$ in the definition of $f_S$ such that the conditions
$\{ 1_h : h \in gS' \} = \{ 1_h : h \in gS \}$ and $\{ 1_t : t \in g S'' \} = \{ 1_t : t \in G \setminus gS \}$ are also
guaranteed. Moreover, since $1 \in S$ we may clearly suppose that $S' \ni 1.$ By (\ref{triviality})
$\af_g(1_{g\m}1_h)=1_{g}1_{gh}$ for $g,h \in G.$ Hence for $g, h_1, \ldots,h_s \in G$ we see that

\begin{align*}
& \af_g(1_{g\m} 1_{h_1} 1_{h_2} \ldots 1_{h_s}) = \af_g((1_{g\m}1_{h_1}) (1_{g\m}1_{h_2})  \ldots (1_{g\m} 1_{h_s})) = \\&
= \af_g(1_{g\m}1_{h_1} )  \af_g(1_{g\m}1_{h_2} )\ldots \af_g(1_{g\m} 1_{h_s})
= (1_{g}1_{gh_1})(1_{g}1_{gh_2}) \ldots (1_{g}1_{gh_s}) = \\&= 1_{g}1_{gh_1} 1_{gh_2} \ldots 1_{gh_s}.
\end{align*}

\noindent Consequently, for $g \in G$ with $g\m \in S$ we obtain

\begin{align*}
& \af_g(f_S) = \af_g(1_{g\m} \prod_{h\in S'}1_h \prod_{t \in S''}(1-1_t)) = 1_{g} \prod_{h\in S'}1_{gh} \prod_{t \in
S''}(1-1_{gt}) = \\&= \prod_{h\in gS'}1_{h} \prod_{t \in gS''}(1-1_{t}) = f_{gS}.
\end{align*}

\noindent Thus  we have

\begin{equation}\label{star4}
\af_g(f_S) = f_{gS} \;\;\;(g\m \in S).
\end{equation}

According to our assumption $St(A_1)=H$ has index $n$ in $A_1.$ Therefore there exist elements $g_1, \ldots, g_n \in G$
with $g_1\m, \ldots g_n\m \in A_1$ such that we have $n$ distinct sets $g_1 A_1, \ldots, g_n A_1.$ It follows from
(\ref{star3})  and (\ref{star4}) that up to a permutation these are exactly the sets $A_1, A_2, \ldots, A_n.$ Hence we may
suppose that $g_i A_1 = A_i\; (i=1, \ldots,n).$ In particular $A_i \neq A_j$ if  $1 \leq i \neq j \leq n.$ Similarly as in the
proof of (\ref{star3}) we easily see that $(f_{A_i})_{jj} =0$ for every $j \neq i$ and consequently $f_{A_i} =
e_{i,i}(1)$ for each $i = 1, \ldots,n.$

 Let $\pi'$ be the elementary partial representation of $G$ determined by the set $A_1.$ Then by (\ref{star2}), $\e'_g =
\pi'(g) \pi'(g\m) = \sum_{A_i \ni g } e_{i,i}(1) =1_g$ for every $g \in G$ and the isomorphism $\af^{\pi'}_g$ has the same
domain as
$\af_g.$ Moreover, the algebra $K^n \equiv diag(K,\ldots,K)$ is generated by the $1_g$'s and for arbitrary $g, h_1, \ldots,
h_s
\in G$ we have

\begin{align*}
& \af^{\pi'}_g(1_{g\m} 1_{h_1} \ldots 1_{h_s}) = \pi'(g) \e'_{g\m} \e'_{h_1} \ldots \e'_{h_s} \pi'(g\m) =
\pi'(g) \e'_{g\m} \pi'(g\m) \e'_{gh_1} \ldots \e'_{gh_s}  = \\&= \e'_{g}  \e'_{gh_1} \ldots \e'_{gh_s}  = \af_g(1_{g\m}
1_{h_1} \ldots 1_{h_s}).
\end{align*}

\noindent Consequently, the partial actions $\af^{\pi'}$ and $\af$ coincide. On the other hand, the  fact that $K^n$ is
generated by the $1_g$'s implies by Proposition~\ref{hom1} that the partial actions $\af$ and $\af^{\pi_{\af}}$ are
equivalent.

It remains to observe that the partial representations $\pi'$ and $\pi_{\af}$ are equivalent. Since $\af^{\pi'} = \af,$
we have by the first part of the proof the isomorphism   $\phi_{\pi'} : K^n \stimes_{\af} G \to  M_n(KH), $ given by  $\sum_{g \in
G} a_g  \de_g \mapsto  \sum_{g \in G} a_g  \pi'(g), (a_g \in \D_g).$ We see that $\phi_{\pi'}(\pi_{\af}(g)) = \phi_{\pi'}(1_g
\de_g) = \phi_{\pi'}(\e'_g \de_g) = \e'_g \pi'(g) = \pi'(g)$ for all $g \in G.$ This yields the equivalence of $\pi'$ and
$\pi_{\af}.$ \fim

\begin{cor}\label{isoBIS}
  For each elementary partial representation  $\pi: G \to M_n(KH)$  the map  $\phi_{\pi} : K^n
\stimes_{\af^{\pi}} G \to  M_n(KH), $ given by  $\sum_{g \in G} a_g  \de_g \mapsto  \sum_{g \in G} a_g  \pi(g),$ is an isomorphism.
\end{cor}

\begin{cor}\label{transitive} The partial action $\af^{\pi}$ which corresponds to an elementary partial
representation  $\pi: G \to M_n(KH)$ acts transitively on the minimal idempotents of $diag(K,\ldots,K),$ i.e. for every
$ 1 \leq i,j \leq n$ there exists an element $g \in G$ such that
$e_{i,i}(1) \in \D_{g\m},  e_{j,j}(1) \in \D_{g}$ and $\af_g(e_{i,i}(1)) = e_{j,j}(1).$
\end{cor}

\p  For arbitrary fixed
$1 \leq i,j \leq n$ there exists  $g \in G$ with $g A_i = A_j.$ Since $g\m \in A_i$ and  $\e_{g\m} = \sum_{\substack{ A_k\ni
g\m}} e_{kk}(1)$ we see that $e_{i,i}(1) \in  \e_{g\m} diag(K,\ldots,K) = \D_{g\m}.$ By (\ref{1star}) and the equality
$\prod_{f \in A_i} \e_f  = e_{ii}(1)$ we have

\begin{align*}
& {\af}^{\pi}_g(e_{i,i}(1)) =  {\af}^{\pi}_g(\prod_{f \in A_i} \e_f) = \pi(g) \prod_{f \in A_i} \e_f \pi(g\m) = \prod_{f
\in A_i} \e_{gf} \pi(g) \pi(g\m) =\\ & = (\prod_{f \in A_j} \e_f) \e_g = \prod_{f \in A_j} \e_f = e_{j,j}(1),
\end{align*}
\noindent as desired. \fim \\

We see that the full matrix algebra  $M_n(KH)$ can be viewed as a crossed product  $K^n \stimes_{\af} G$ in many ways. In particular,
we have the following.

\begin{cor}\label{particular} For each positive integer $n$ and an arbitrary group $G$ of order $n+1$ there is a partial
action
$\af$ of $G$ on $K^n$ such that  $M_n(K) \cong K^n \stimes_{\af} G.$
\end{cor}

\p Pick an element $1 \neq g \in G$ and take $A = G \setminus \{g\}.$ Then $St(A) = 1$ and $A$ determines an elementary partial
representation $\pi: G \to M_n(K),$ which gives rise to the isomorphism $\phi_{\pi} : K^n \stimes_{\af^{\pi}} G \to  M_n(K). $ \fim
\\

A crossed product structure on  $M_n(KH)$ gives a grading of this $K$-algebra by $G.$ It is easily seen that this grading
is elementary. More precisely, we recall that a grading on the $K$-algebra $M_n(K)$ by a group $G$ is called {\it elementary}
if each elementary matrix
$e_{i,j}(1)$ is homogeneous. The elementary gradings are essential for the description of gradings on matrix algebras (see
\cite{BSZ}). It is known that each elementary  grading on
$M_n(K)$ by a group
$G$ is determined by an
$n$-tuple $(g_1=1, g_2,
\ldots, g_n)$ of non-necessarily distinct elements of $G$ in such a way that the homogeneous degree $ deg( e_{i,j}(1))$ of $
e_{i,j}(1)$ is
$g_i\m g_j.$ Conversely, in this manner each $n$-tuple $(g_1=1, g_2, \ldots, g_n)$ determines an elementary grading on
$M_n(K).$ It turns out that if in  $(g_1=1, g_2, \ldots, g_n)$ the $g_i$'s are pairwise distinct and
 $St(\{g_1, g_2, \ldots, g_n\})=1,$ then the corresponding elementary grading of $M_n(K)$  necessarily comes from a
crossed product structure  $K^n \stimes_{\af^{\pi}} G \cong M_n(K).$ More generally, we  say that a grading of the
$K$-algebra $M_n(KH)$ by a group
$G$ is {\it elementary} if for each $h \in H, i,j
\in \{1, \dots, n\}$ the elementary matrix $e_{i,j}(h)$ is homogeneous. Given a subset $1 \in A$ of $G$ with $St(A)=H$ one
defines an elementary grading of $M_n(KH)$  by the equality $deg(e_{i,j}(h)) = g_i\m h g_j\; (h \in H, i,j \in \{1, \dots,
n\}),$ where  $A = Hg_1 \cup Hg_2 \cup \ldots \cup Hg_n, \; g_1 = 1$  and $Hg_i \neq Hg_j$ for $1 \leq i \neq j \leq n.$
Changing the
$g_i$'s to $g_i\m$'s in the definition of the elementary partial representations, we easily obtain the following.

\begin{cor}\label{gradings}  For the
elementary grading
of $M_n(KH)$ by  a group $G,$ determined by a  subset $ A  \subseteq G$ with $A \ni 1,$
$H=St(A),$ and for the
ele\-men\-ta\-ry partial representation\\ $\pi: G \to M_n(KH),$    constructed from $A_1=A, A_2  = g_2\m A, \ldots,
A_n = g_n\m
A,$ the map $\phi_{\pi} : K^n\stimes_{\af^{\pi}} G \to  M_n(KH) $ is an isomorphism  of graded $K$-algebras.

\end{cor}

\end{section}

\begin{center}
 \large{Acknowledgments}
\end{center}

While this work was carried out the first author  visited several times the Federal University of Santa Catarina (Brazil) in 2001-2003. He expresses his appreciation to the Department of Mathematics of that University for its warm hospitality and financial support during the visits.


\begin{thebibliography}{99}

\bibitem{Abadie} F.\ Abadie, Sobre a\c c\~oes parcias, fibrados de
Fell e grup\'oides, PhD Thesis, Universidade de S\~ao Paulo, 1999.

\bibitem{AbadieTwo} F. Abadie, Enveloping Actions and Takai Duality
for Partial Actions, {\it J.~Funct.~Analysis}, to appear,
[arXiv:math.OA/0007109].

\bibitem{BSZ} Yu.\ Bahturin, S.\ K.\ Sehgal, M.\ V.\  Zaicev, Group gradings on associative algebras, {\it J. Algebra},
{\bf 241} (2), (2001), 677-698.


\bibitem{CK} J.\ Cuntz, W.\ Krieger,  A Class of $C^*$-Algebras and Topological
Markov Chains, {\it Inventiones Math.}  {\bf 56} (1980), 251--268.

\bibitem{DEP} M.\ Dokuchaev, R.\ Exel, P.\ Piccione, Partial representations
and partial group algebras, {\it J. Algebra}, {\bf 226} (1), (2000), 505-532.

\bibitem{DZh0} M.\ Dokuchaev, N. Zhukavets, On finite degree    partial representations of groups,  {\it J. Algebra}, to appear.

\bibitem{DZh} M.\ Dokuchaev, N. Zhukavets, On irreducible  partial representations of groups,
{\it Comptes Rendus Math. Rep.  Acad.  Sci.  Canada}, {\bf 24} (2), 2002, 85-90.

\bibitem{E-1} R.\ Exel, Circle actions on $C^*$-algebras, partial automorphisms and
generalized Pimsner-Voiculescu exact sequences,  {\it J. Funct. Anal.} {\bf 122} (3),
(1994), 361 - 401.

\bibitem{E0} R.\ Exel, Twisted partial actions: a classification of regular
$C^*$-algebraic bundles, {\it Proc.\ London Math. Soc.} {\bf 74} (3),  (1997), 417 - 443.

\bibitem{E1} R.\ Exel,  Partial Actions of Groups and Actions of
Semigroups,  {\it Proc.\ Am.\ Math.\ Soc.} {\bf 126} (12),  (1998), 3481--3494.

\bibitem{E2} R.\ Exel,  Amenability for Fell Bundles, {\it J.\ Reine Angew.\ Math.}
{\bf 492} (1997), 41--73.

%\bibitem{E3} R.\ Exel,  Partial Representations and Amenable Fell Bundles
%over Free Groups,  {\it Pacific J.\ Math.} {\bf 192} (1), (2000), 39-63.

\bibitem{EL} R.\ Exel, M.\ Laca,  Cuntz--Krieger Algebras for Infinite Matrices,
{\it J.\ Reine Angew.\ Math.} {\bf 512} (1999), 119-172.

\bibitem{FellDoran} J.\ M.\ G.\ Fell, P.\ S.\ Doran, {\it Representations of *-Algebras, Locally
Compact Groups and Banach *-Algebraic Bundles I}, Academic Press, 1988, Pure and Applied
Math. {\bf 125}.

\bibitem{Fillmore} P.\ A.\ Fillmore, {\it A User's Guide to Operator Algebras,}
Willey - Interscience, 1996.



\bibitem{Mc} K.\ McClanahan, K-theory for partial crossed products by discrete groups,
{\it J. Funct. Anal.} {\bf 130} (1995), 77 - 117.

 \bibitem{Passman} D.\ S.\ Passman, {\em The Algebraic Structure of Group Rings},
 Interscience, New York, 1977.

\bibitem{QR} J.\ C.\ Quigg, I.\ Raeburn,  Characterizations of Crossed
Products by Partial Actions, {\it J.\ Operator Theory} {\bf 37} (1997), 311--340.


\bibitem{Rowen} L.\ H.\ Rowen, {\em Ring theory - Student edition}, Academic Press, 1991.

 \bibitem{Sehgal} S.\ K.\ Sehgal, {\em Units in Integral Group Rings}, Longman
 Scientific \& Technical Press, Harlow, 1993.

\end{thebibliography}
\end{document}